\documentclass[12pt]{amsart}
\usepackage{amssymb}
\newtheorem{theorem}{Theorem}[section]

\newtheorem{proposition}[theorem]{Proposition}

\newtheorem{lemma}[theorem]{Lemma}
\theoremstyle{definition}
\newtheorem{definition}[theorem]{Definition}
\theoremstyle{remark}
\newtheorem{remark}[theorem]{Remark}

\newcommand\A{\mathcal{A}}
\newcommand\M{\mathcal{M}}

\renewcommand{\L}{\mathcal{L}}
\renewcommand{\O}{\mathcal{O}}

\newcommand{\U}{\mathcal{U}}

\newcommand{\N}{\mathbb{N}}
\newcommand{\R}{\mathbb{R}}
\newcommand{\RR}{\mathcal{R}}
\newcommand{\C}{\mathbb{C}}
\newcommand{\cC}{\mathcal{C}}

\newcommand{\Z}{\mathbb{Z}}
\newcommand{\Q}{\mathbb{Q}}
\renewcommand{\P}{\mathbb{P}}
\newcommand\lie[1]{\mathfrak{#1}}
\renewcommand{\k}{\lie{k}}

\newcommand{\g}{\lie{g}}
\newcommand{\p}{\lie{p}}

\renewcommand{\t}{\lie{t}}
\newcommand{\Alc}{\lie{A}}

\newcommand{\on}{\operatorname}

\newcommand{\Aut}{ \on{Aut} } 
\newcommand{\Bun}{ \on{Bun} } 
\newcommand{\ParBun}{ \on{ParBun} } 
\newcommand{\ParVect}{ \on{ParVect} } 
 
\newcommand{\Ad}{ \on{Ad} }

\newcommand{\Hol}{ \on{Hol} } 

\newcommand{\Hom}{ \on{Hom}}

\newcommand{\Spec}{  \on{Spec}}

\newcommand\dirac{/\kern-1.2ex\partial} % Dirac operator
\newcommand\qu{/\kern-.7ex/} % Categorical quotients
 
\newcommand{\labell}\label

\renewcommand{\d}{{\mbox{d}}}
\newcommand{\ol}{\overline}

\newcommand\Phinv{\Phi^{-1}}

\newcommand\ddt{\frac{\d }{\d t}}

\newcommand\eps{\epsilon}

\newcommand{\Del}{\Delta}

\newcommand{\ti}{\tilde}

\newcommand\SL{\on{SL}}

\newcommand\mE{\mathcal{E}}

\newcommand\pardeg{\on{pardeg}}

\renewcommand{\ss}{\on{ss}}
\newcommand\Tr{\on{Tr}}

\newcommand\Gl{\on{Gl}}
\newcommand\Gr{\on{Gr}}
\newcommand\Map{\on{Map}}
\newcommand\rank{\on{rank}}

\begin{document}

\title[Parabolic bundles and products of conjugacy classes]{Parabolic
bundles, products of conjugacy classes, and Gromov-Witten invariants}

\author{C. Teleman}
\address{DPMMS, CMS, Wilberforce Road, Cambridge, CB3 0WB, UK}
\email{teleman@dpmms.camac.uk}

\author{C. Woodward}
\address{Mathematics-Hill Center, Rutgers University,
110 Frelinghuysen Road, Piscataway NJ 08854-8019, USA}
\email{ctw@math.rutgers.edu}

\begin{abstract}
The set of conjugacy classes appearing in a product of conjugacy
classes in a compact, $1$-connected Lie group $K$ \label{Ksec} can be
identified with a convex polytope in the Weyl alcove \cite{me:lo}.  In
this paper we identify linear inequalities defining this polytope.
Each inequality corresponds to a non-vanishing Gromov-Witten invariant
for a generalized flag variety $G/P$, where $G$ is the
complexification of $K$ and $P$ is a maximal parabolic subgroup.  This
generalizes the results for $SU(n)$ of Agnihotri and the second author
\cite{ag:ei} and Belkale \cite{bl:ip} on the eigenvalues of a product
of unitary matrices and quantum cohomology of Grassmannians.
\end{abstract}



\maketitle

\section{Introduction}

An old problem which goes back to Weyl is to determine the possible
eigenvalues of a sum of traceless Hermitian matrices.  According to a
result of Klyachko \cite{kl:st}, see also \cite{bl:ip,kt:ho2}, there
is a finite set of homogeneous linear inequalities on the eigenvalues,
each of which corresponds to a non-vanishing structure coefficient in
the Schubert calculus of a Grassmannian.  The same inequalities turn
out to determine the non-vanishing of the Littlewood-Richardson
numbers \cite{kt:sa}.  Berenstein-Sjamaar \cite{be:coa} and
Leeb-Millson \cite{mi:co} generalize this result to arbitrary type as
follows.  Let $\k$ be the Lie algebra of $K$, and $T$ a maximal torus
\label{Tsec} with Lie algebra $\t$.  The set of coadjoint orbits in
$\k^*$ is parametrized by a Weyl chamber $\t^*_+$ in the fixed point
set $\t^*$ of the action of $T$ on $\k^*$.  For any $\mu \in \t^*_+$,
we denote by $\O_\mu$ the corresponding coadjoint orbit $\O_\mu = K
\cdot \mu$.  For any $\mu_1,\ldots, \mu_{b-1} \in \t^*_+$ 
the sum 
$$ \O_{\mu_1} + \ldots + \O_{\mu_{b-1}} = \bigcup_{\mu_b} \O_{\mu_b}
$$
for some set of $\mu_b$ in $\t^*_+$.  Which $\mu_b$'s occur is
determined by a finite number of linear inequalities, each of which
corresponds to a non-vanishing structure coefficient in the Schubert
calculus for $G/P$.  There are similar results for other symmetric
spaces.

Biswas \cite{bi:r2}, Agnihotri-Woodward \cite{ag:ei} and Belkale
\cite{bl:ip} generalize Klyachko's result to eigenvalues of products
of special unitary matrices.  For any special unitary matrix the
logarithms $\lambda_1,\ldots,\lambda_r$ of the eigenvalues may be
chosen so that
$$\lambda_1 + \ldots + \lambda_r = 0, \ \ 
\lambda_1 \ge \lambda_2 \ge \ldots \ge \lambda_r \ge \lambda_r - 2\pi
.$$ 
There are a finite set of linear inequalities on the $\lambda_i$'s for
a product, each of which corresponds to a non-vanishing structure
coefficient of the {\em quantum} Schubert calculus of a Grassmannian.

In this paper we solve the multiplicative problem for arbitrary type;
this includes as a special case some of the results of
Berenstein-Sjamaar and Leeb-Millson.  Let
\label{chambsec} $\alpha_0 \in \t^*$ \label{alphasec} denote the 
highest root.  The set of conjugacy classes in $K$ is parametrized by
the Weyl alcove
\begin{equation} \label{alcove}
\Alc = \{ \xi \in \t_+, \ \alpha_0(\xi) \leq 1 \} .\end{equation}
For any $\mu \in \Alc$, we denote by $\mathcal{C}_\mu$ the conjugacy
class of $\exp(\mu)$.  For $\mu_1,\ldots,\mu_{b-1}$, the product $
\cC_{\mu_1} \cdot \cC_{\mu_2} \ldots \cdot \cC_{\mu_{b-1}}$ is
invariant under conjugation; we wish to identify which conjugacy
classes $\cC_{\mu_b}$ appear in
$$ \cC_{\mu_1} \cdot \cC_{\mu_2} \ldots \cdot \cC_{\mu_{b-1}}
= \bigcup_{\mu_b} \cC_{\mu_b} .$$
More symmetrically, define
$$ \Del_b = \{ (\mu_1,\ldots,\mu_b) \in \Alc^b \ | \
\mathcal{C}_{\mu_1} \ \cdot\ldots \cdot \mathcal{C}_{\mu_b} \ni e \}
$$
where $e$ is the group unit.  By \cite[Corollary 4.13]{me:lo},
$\Del_b$ is a convex polytope of maximal dimension in $\Alc^b$.  We
wish to find the defining inequalities for $\Del_b$.

The polytope $\Del_b$ can also be described as the possible holonomies
of flat $K$-bundles on the punctured two-sphere.  Let $x_1,\ldots,x_b$
be distinct points on an oriented surface $X$
\label{Xsec}.  For any markings $\mu_1,\ldots,\mu_b$ there exists a
symplectic stratified space $\RR_K(X;\mu_1,\ldots,\mu_b)$ whose points
are the isomorphism classes of flat $K$-bundles on $X \backslash \{
x_1,\ldots,x_b \}$ with holonomy around $x_i$ in $\cC_{\mu_i}$.
Equivalently, $\RR_K(X;\mu_1,\ldots,\mu_b)$ is the moduli space of
representations of the fundamental group mapping a small loop around
$x_i$ to $\cC_{\mu_i}$.  In the case $X$ has positive genus
$\RR_K(X;\mu_1,\ldots,\mu_b)$ is always non-empty. In the genus zero
case we have
$$ \Delta_b = \{(\mu_1,\ldots,\mu_b), \ \ \RR_K(X;\mu_1,\ldots,\mu_b)
\neq \emptyset \} .$$

A final interpretation of the problem involves the space of conformal
blocks, or equivalently, fusion products of representations of affine
Lie algebras.  Fix a complex structure on $X$.  The space
$\RR_K(X;\mu_1,\ldots,\mu_b)$ may be identified with the moduli space
of semistable parabolic bundles on $X$, by a theorem of Mehta-Seshadri
\cite{ms:pb}, Bhosle-Ramanathan \cite{bh:pa}, and the discussion in
Section 4.  Let
$$ \t_\Q = \Lambda^* \otimes_\Z \Q $$ 
be the set of rational points in the Cartan, and suppose
$(\mu_1,\ldots,\mu_b) \in \t_Q$. Then there exists $n \in \N$ such
that $n\mu_i$ are all dominant weights.  The basic line bundle over
$\RR_K(X;\mu_1,\ldots,\mu_b)$ determines a projective embedding with
Hilbert polynomial given by the dimension of the space of genus zero
conformal blocks ${\mathcal H}(X;kn\mu_1,\ldots,kn\mu_b;k)$ at level
$k$ with markings $kn\mu_1,\ldots,kn\mu_b$
\cite{pa:co},\cite{la:li},\cite[Section 8]{te:qu}.
Thus 
$$\Delta_b \cap \t_{\Q} = \{ (\mu_1,\ldots,\mu_b), \ \ \exists k
\ \text{such that}\ 
\mathcal{H}(X;kn\mu_1,\ldots,kn\mu_b;k) \neq \{ 0 \} \} .$$

Our description of the inequalities for $\Delta_b$ involves the small
quantum cohomology $ QH^*(G/P)$, a deformation of the ordinary
cohomology ring defined by including contributions from higher degree
rational curves in $G/P$.  For simplicity, we discuss only the case
that $P$ is maximal.  Recall that the Schubert basis for $H^*(G/P)$ is
given by the classes of closures of orbits of a Borel subgroup on
$G/P$.  Let $B$ be the standard Borel subgroup whose Lie algebra
contains the positive root spaces.  Let $P \subset G$ denote a
parabolic subgroup, corresponding to a subset $\Pi_P$ of the simple
roots $\Pi$.  Let $W_P \subset W$ the subgroup of $W$ generated by
simple reflections for roots $\alpha \in \Pi_P$. \label{Bsec} For any
$w \in W/W_P$, the {\em Schubert variety}
\label{Ysec}
$$Y_w = \ol{BwP/P} \subset G/P$$
is a normal subvariety of $G/P$.  The homology classes $ [Y_w]$ form a
basis for the homology $H^*(G/P)$; in this paper we use rational
coefficients.  Let $w_o \in W$ be the long element in the Weyl group.
The class of $Y^w := Y_{w_ow}$ is Poincar\'e dual to $[Y_w]$.  Its
degree is $\deg [Y^w] = 2l_P(w)$, where $l_P(w)$ is the minimal length
of a representative of $w$ in $W$.  Now let $q$ be a formal variable.
As a $\Q[q]$-module $QH^*(G/P)$ is freely generated by $H^*(G/P).$ Fix
$X = \P^1$ and choose distinct points $x_1,\ldots,x_b \in
X$. \label{ndsec} For any holomorphic map $\varphi: \, X \to G/P$ the
degree of $\varphi$ is
$$\deg(\varphi) := \varphi_* [X] \in H_2(G/P)
\cong \Z .$$ 
Let $g_i Y_{w_i}, i = 1,\ldots,b$ be general translates of the
Schubert varieties $ Y_{w_i}$.  Let $n_d(w_1,\ldots,w_b)$ be the
number of holomorphic maps $\varphi: \ X \to G/P$ of degree $d$ such
that $\varphi(x_i) \in g_iY_{w_i}$, if this number is finite, and zero
otherwise.  Define
$$ [Y_{w_1}] \star \ldots \star [Y_{w_{b-1}}] = \sum_{d \in \N, \ w_b
\in W/W_P} n_d(w_1,\ldots,w_b) q^d [Y^{w_b}].
$$
The resulting product is commutative, associative, and independent of
the choice of $x_1,\ldots,x_b$ and general $g_1,\ldots,g_b$
\cite{fu:st,fu:qp}.  These Gromov-Witten invariants of $G/P$ (as
opposed to the invariants that appear in the large quantum cohomology)
are computable in practice using formulas of D. Peterson \cite{pe:mi},
whose proofs are given in \cite{fu:qp} and \cite{wo:pe}.  An example,
for the case $G_2$, is given at the end of the paper.

For any maximal parabolic subgroup $P$, let $\omega_P$ denote the
fundamental weight that is invariant under $W_P$.  Our main result is

\begin{theorem} \label{main}  The polytope $\Del_b$ is the set 
of points $(\mu_1,\ldots,\mu_b) \in \Alc^b$ satisfying
$$ \sum_{i=1}^b (w_i \omega_P,\mu_i ) \leq d $$
for all maximal parabolic subgroups $P \subset G$ and all
$w_1,\ldots,w_b \in W/W_{P}$ and non-negative integers $d$ such that
the Gromov-Witten invariant $n_d(w_1,\ldots,w_b) = 1$.
\end{theorem}
A connection between this problem and the Hofer metric on
symplectomorphism groups is discussed by Entov \cite{en:ka}.

There are several remaining open questions.  We do not know which
inequalities are independent.  Also, the quantum generalization of the
saturation conjecture \cite{kt:sa}: are the inequalities necessary and
sufficient conditions for the non-vanishing of the fusion
coefficients, at least in the simply laced case?  There are similar
polytopes for products of conjugacy classes in disconnected groups.
These might be related to the twisted quantum cohomology (Floer
cohomology for symplectomorphisms not isotopic to the identity.)

\subsection{Index of notation}

\begin{tabbing}
\indent \= xxxxxxxxxxxx \= xxxxxxxxxxxxxxxxxxxxxxxxxxxxxxxxxxxxxxxxxxxxxxxxxx \= xxxxxx \kill % sample line
\> $K,G$ \> simple $1$-connected compact group, resp. complexification\> p. \pageref{Ksec} \\
\> $T,\t$ \> maximal torus, resp. Cartan subalgebra\> p. \pageref{Tsec} \\
\> $\alpha_0,\t_+,\Alc$  \> highest root, resp. positive chamber, resp. alcove\> p. \pageref{alphasec} \\
\> $\mu_j$ \> marking in $\t_+$ with $\alpha_0(\mu_j) < 1$ \> p. \pageref{parabsec} \\
\> $B,P$\> Borel, resp. standard parabolic subgroup \>  p. \pageref{Bsec} \\
\> $(P_1,P_2)$ \> relative position of parabolics\> p. \pageref{relsec} \\
\> $X$ \> smooth curve $/ \C$\> p. \pageref{Xsec} \\
\> $Y_w, C_w$ \> Schubert variety, resp. Schubert cell \> p. \pageref{Ysec},
p. \pageref{Csec} \\
\> $n_d(w_1,\ldots,w_d)$ \> Gromov-Witten invariant\> p. \pageref{ndsec} \\
\> $\mE \to X$ \> holomorphic principal $G$-bundle\> p. \pageref{Esec} \\
\> $\varphi_j \in \mE_{x_j}/P_j$ \> parabolic reduction of $\mE$ at $x_j$\> p. \pageref{parabsec} \\
\> $\sigma: \ X \to \mE/P$ \> parabolic reduction\> p. \pageref{sigsec} \\
\> $\pi: \ \ti{X} \to X$ \> ramified cover\> p. \pageref{pisec} \\
\> $\ti{x}_j$, $\ti{U}_j$  \> ramification point of $\pi$,
resp. neighborhood
of $\ti{x}_j$\> p. \pageref{pisec}  \\
\> $L,U$\> Levi, resp. unipotent subgroup\> p. \pageref{levisec} \\
\> $r: \ P \to L$ \> projection to $L$ \> p. \pageref{rsec}\\
\> $\iota: \ L \to G $ \> inclusion of $L$ \>  p. \pageref{rsec} \\
\> $ \Lambda_P^*$ \> weights of characters of $P$\> p. \pageref{charsec} \\
\> $\sigma_{\mE},\mu_{\mE}$ \> canonical reduction, slope \>  p. \pageref{canred} \\
\> $\U_G(X)$ \> universal space for $G$-bundles\> p. \pageref{UGsec} \\
\> $\M_G(X;x;\mu)$ \> moduli space of parabolic semistable $G$-bundles\> p. \pageref{modsec} \\
\> $\RR_K(X;\mu)$ \>moduli space of flat $K$-bundles with fixed
holonomy \>  p. \pageref{repsec}   \\
\> $A,A_\infty$\> a connection, resp. its Yang-Mills limit \>  p. \pageref{limsec}  \\
\end{tabbing}

\section{Parabolic $G$-bundles}
 \label{Esec}

In this section we develop the general theory of parabolic
$G$-bundles: equivalence with equivariant bundles for a finite group,
canonical reductions, and coarse moduli spaces.  Unfortunately, we
could not understand the arguments in Bhosle-Ramanathan \cite{bh:pa}
which covers similar material so we chose to employ a different
approach, basically switching the order of embedding $G$ in $GL(n)$
and applying the equivalence with equivariant bundles.  A different
approach which is less useful for our purposes but works in any
dimension is given by Balaji, Biswas, and Nagaraj \cite{ba:pr}.

\subsection{Definitions}

Let $X$ be a complex manifold.  A principal $G$-bundle over $X$ is a
complex manifold $\mE \to X$ with a right action of $G$ that is locally
trivial.  That is, any point in $X$ is contained in a neighborhood $U$
such that $\mE|_U$ is $G$-equivariantly biholomorphic to $U \times G$.
For $X$ a scheme, principal $G$-bundles over $X$ are required to be
locally trivial in the \'etale topology.  By the theorem of Drinfeld
and Simpson \cite{dr:bs} (another proof is given in \cite{te:bo}) any
principal $G$-bundle over the product $X \times S$ of a smooth curve
$X$ with a scheme $S$ is trivial locally in the product of the Zariski
topology for $X$ and the \'etale topology in $S$.  The results of this
section are mostly valid in both the analytic and algebraic
categories.

The following definition of parabolic vector bundle is slightly more
general than the original one given by Mehta and Seshadri.  Let $X$ be
a curve with distinct marked points $x_1,\ldots,x_b$, and $E \to X$ a
holomorphic vector bundle of rank $r$.  A {\em parabolic structure}
for $E$ at $x_i$ is a partial flag
$$ E_{x_i}^1 \subset E_{x_i}^2 \subset \ldots \subset E_{x_i}^{l_i} =
E_{x_i} $$
together with a set of {\em markings}
$$ \mu_{i,1} \ge \mu_{i,2} \ge \ldots \ge \mu_{i,r}, \ \ \mu_{i,1} -
\mu_{i,r} < 1 $$
corresponding to the type of the partial flag.  That is, for all $j =
1,\ldots,r$,
$$\# \{ \mu_{i,k}, \ k \leq j \} = \# \{ \dim(E_{x_i}^k)\leq j \} .$$
Note that Mehta-Seshadri require that the markings be non-negative, so
that the $\mu_i= (\mu_{i,1}, \ldots, \mu_{i,r})$ lie in a fundamental
domain for the affine Weyl group of $\lie{gl}(r)$.  

A parabolic vector bundle over a pointed curve $(X;x_1,\ldots,x_b)$ is
a holomorphic vector bundle $E \to X$ together with parabolic
structures $(E_{x_i}^\bullet, \mu_i)$ at the points $x_i, \ i =
1,\ldots,b$.  Usually we drop the parabolic structures from the
notation.

We define a parabolic $\SL(r)$-vector bundle to be a parabolic vector
bundle $E \to X$ with degree zero and
$$ \sum_{j=1}^r \mu_{i,j} = 0, \ \ i = 1,\ldots, b.$$
Hence the markings $\mu_i$ lie in the Weyl alcove \eqref{alcove} for
the Lie algebra $\lie{sl}(r)$.  There is an equivalent definition of
parabolic structure in terms of the bundle $\mE$ of frames for $E$.
For any $\mu \in \t_+$, there is a unique standard parabolic subgroup
$P \subset \SL(r)$ such that the $W_P$ is the stabilizer of $\mu$ in
$W$.  We say that $P$ is the standard parabolic subgroup {\em
corresponding to } $\mu$.  Let $P_{i}$ denote the standard parabolic
subgroup corresponding to the marking $\mu_i$.  The data of the
filtration of $E_{x_i}$ is equivalent to a reduction of $\mE_{x_i}$ to
the parabolic subgroup $P_{i}$.  Explicitly, let $\varphi_i$ denote
the set of frames $\{v_1,\ldots,v_r\}$ for $E_{x_i}$, such that $v_l
\in E_{x_i}^j$ for $l \leq \dim(E_{x_i}^j)$.  Then $P_{i}$ acts
transitively on $\varphi_i$, that is, $\varphi_i$ is a reduction of
$\mE_{x_i}$ to structure group $P_{i}$.

Let $\mE \to X$ be a principal $G$-bundle.

\begin{definition}  A {\em parabolic structure} for
$\mE$ at $x_i$ consists of
\begin{enumerate}   \label{parabsec}
\item a {\em {marking}} $\mu_i \in \Alc$ with $\alpha_0(\mu_i) <1$;
\item a reduction $\varphi_i \in \mE_{x_i}/P_{i}$, where $P_{i}$ is
the standard parabolic subgroup corresponding to $\mu_i$.
\end{enumerate}
A {\em parabolic bundle} on $(X;x_1,\ldots,x_b)$ is a bundle $\mE$ with
parabolic structures at $x_1,\ldots,x_b$.  A {\em family of parabolic
bundles} parametrized by a complex manifold $S$ is a principal
$G$-bundle over $X \times S$ with sections of $\mE/P_i$ over $\{ x_i \}
\times S$ and markings $\mu_i$.  A morphism of bundles $\mE_1 \to \mE_2$
defines a morphism of parabolic bundles if the bundles have the same
markings and parabolic reductions for $\mE_1$ are mapped to parabolic
reductions for $\mE_2$.
\end{definition} 

We remark that one can replace the condition $\alpha_0(\mu_i) < 1$
with $\alpha_0(\mu_i) \leq 1$, by working with torsors (non-abelian
cohomology classes) for a group sheaf which is locally a standard
parabolic subgroup of the loop group; see \cite{te:bo} for
definitions.  However, we do not know any intrinsic formulation of the
semistability condition in this language.  In the case $G= \SL(r)$,
all parabolic subgroups of the loop group are conjugated to subgroups
of $G[[z]]$ by outer automorphisms. That is why this case does not
need to be considered for moduli spaces of vector bundles.

The parabolic degree of a parabolic vector bundle $E$ is defined by
$$ \pardeg(E) := \deg(E) + \sum_{i=1}^b \sum_{j = 1}^{r} \mu_{i,j} .$$
Here $\deg(E)$ denotes the first Chern class $c_1(E) \in H_2(X) \cong
\Z$.  The parabolic slope of $E$ is
$$\mu(E) := \pardeg(E)/\rank(E) .$$
The parabolic structure on $E$ induces a parabolic structure on any
holomorphic subbundle $F$.  Define a flag in $F_{x_i}$ by removing
repeating terms from the sequence 
$$F_{x_i} \cap E^1_{x_i} \subseteq F_{x_i} \cap E^2_{x_i} \subseteq
\ldots \subseteq F_{x_i} \cap E^{l_i}_{x_i}.$$
Define markings $\nu_{i}$ by $\nu_{i,j} = \mu_{i,k}$ where
$k$ is the smallest integer such that $F^j_{x_i} \subset E^k_{x_i}$.
The parabolic bundle $E$ is {\em semistable} if and only if the
inequality
$$ \mu(F) \leq \mu(E) $$
holds for all subbundles $F \subset E$.

In order to generalize these definitions to arbitrary type we give a
definition of ordinary semistability using the frame bundle $\mE$ for
$E$.  Let $P_k$ denote the standard maximal parabolic subgroup of
$GL(r)$ stabilizing a subspace of dimension $k$.  Let $\sigma: \ X \to
\mE/P_k$ be the parabolic reduction with $\sigma(x)$ equal to the set
of frames for $E_x$ whose first $k$ elements are in $F_x$.  Let
$\mE(\omega_k)$ denote the line bundle $\mE \times_{P_k} \C_{\omega_k}
$ where $\C_{\omega_k}$ is the weight space for the $k$-th fundamental
weight $\omega_k$ of $GL(r)$.  A little yoga with the definition of
Chern classes shows that
$$ \deg(F) = \deg(\sigma^* \mE(\omega_k)) .$$
If $E$ is an $\SL$-vector bundle, then $E$ is ordinary semistable if
and only if
$$ \deg(\sigma^* \mE(\omega_k)) \leq 0, \ \ \forall \sigma: \ X \to
\mE/P_k $$
for all $k = 1,\ldots,r-1$.  The definition of the marking $\nu$ can
be rephrased in terms of the Schubert cell decomposition of the
Grassmannian.  Let $V$ be a vector space of dimension $r$, and
$$ V^\bullet = \{ V^1
\subset V^2 \subset \ldots V^r = V \} $$ 
a complete flag in $V$.  Let $\Gr(k,V)$ denote the Grassmannian of
$k$-planes in $V$.  For each sequence of integers $J = \{ j_1 < \ldots
<j_k \}$ the Schubert variety corresponding to $I$ is
$$ Y_I = \{ U \subset V, \ \ \dim(U) \cap V^{j_m} \ge m, \ \ m = 1,\ldots, k
.\}
$$
Let $C_I$ denote the interior of $Y_I$, that is, $C_I = Y_I \backslash
\cup Y_J$, where the union is over $Y_J$ contained in $Y_I$.  We say
that $U$ is in {\em relative position} $I$ to $V^\bullet$ if $U$ lies
in $C_I$.  Now let $E \to X$ be a parabolic vector bundle, and $F
\subset E$ a holomorphic sub-bundle.  The marking for $F_{x_i}$ is
$(\mu_{j}, j \in J_i)$ where $J_i$ is the relative position of
$F_{x_i}$ and $E_{x_i}^\bullet$.  We can write this in the language of
principal bundles as follows.  The quotient $\mE_{x_i}/P_{i}$ is
isomorphic to the flag variety for $E_{x_i}$ of type corresponding to
$\mu_i$, so the flag $E_{x_i}^\bullet$ defines a point $\varphi_i \in
\mE_{x_i}/P_{i}$.  The quotient $\mE_{x_i}/P$ is isomorphic to the
Grassmannian $\Gr(k,E_{x_i})$, and any subspace $U_{x_i} \subset V$
defines a point $\sigma(x_i) \in \mE_{x_i}/P_k$.  The quotient of Weyl
groups $W/W_{P_k}$ maps bijectively to the set of elements of size $k$
in $\{1,\ldots,r\}$, by
$$[w] \mapsto I([w]) := \{ w(k+1),\ldots, w(r)\} .$$
We say that $\sigma(x_i)$ is in relative position $[w_i]$ to
$\varphi_i$ if $U_{x_i}$ is in relative position $I([w_i])$ to
$E_{x_i}^\bullet$.  Hence $E$ is parabolic semistable if and only if
$$ \deg(\varphi^* \mE(\omega_k)) + \sum_{i = 1}^b (w_i \omega_k,\mu_i)
\leq 0 $$
for all $k=1,\ldots,r-1$ and reductions $\varphi: \ X \to \mE/P_k$,
where $[w_i] \in W_{P_i} \backslash W/W_P$ is the relative position of
$\sigma(x_i)$ and $\varphi_i$, and $w_i$ is any representative of
$[w_i]$ in $W$.

For arbitrary simple $G$ and parabolic subgroups
\label{relsec} $P_1' = \Ad(g_1)P_1, P_2' = \Ad(g_2)P_2 \subset G$
given as conjugates of standard parabolics $P_1$ and $P_2$, define the
{\em relative position} $(P_1', P_2') \in W_{P_1} \backslash W /
W_{P_2}$ to be the image of $(g_1,g_2)$ under the map
$$ G \times G \to G \backslash (G \times G) / P_1 \times P_2 \cong
{P}_{1} \backslash G / {P}_{2} \cong W_{{P}_1} \backslash W /
W_{{P}_2} .$$
Note that $ (P_2',P_1') = (P_1',P_2')^{-1}$ and $(P',P') = [1] $ for
any parabolic subgroup $P'$.  

\begin{definition} A parabolic principal $G$-bundle 
$(\mE;\mu_1,\ldots,\mu_b;\varphi_1,\ldots,\varphi_b)$ is {\em stable}
(resp. {\em semistable}) if for any maximal parabolic subgroup $P$ and
reduction \label{sigsec} $\sigma: \ X \to \mE/P$ we have
\begin{equation} \label{ineq}
 \deg(\sigma^*\mE(\omega_P)) + \sum_{i=1}^b (w_i \omega_P, \mu_i) < 0 \
\text{ (resp. } \ \leq 0 \text{ \ ) } .\end{equation}
where $w_i \in W_{P_i} \backslash W/W_{P}$ is the relative position of
$\sigma(x_i)$ and $\varphi_i$.
\end{definition}
\noindent By $(w_i \omega_P,\mu_i)$ we mean $(\ti{w_i}
\omega_P,\mu_i)$, independent of the choice of representative
$\ti{w_i}$ of $w_i$.  We call the left-hand-side of \eqref{ineq} the
{\em parabolic degree} of $\sigma$.

If $G = \SL(r)$ and $\mE$ is a parabolic principal $G$-bundle, then
the parabolic structure induces on the associated vector bundle $E$
the structure of a principal $\SL(r)$-vector bundle.  For any smooth
curve $X$ with marked points $x_1,\ldots,x_b$, let
$\ParVect_0(X;x_1,\ldots,x_b)$ denote the functor which assigns to any
complex manifold $S$ the set of isomorphism classes of families of
parabolic $\SL(r)$-vector bundles on $(X;x_1,\ldots,x_b)$ parametrized
by $S$.  Let $\ParBun(X;x_1,\ldots,x_b;G)$ denote the functor that
assigns to any complex manifold $S$ the set of isomorphism classes of
families of parabolic principal $G$-bundles on $(X;x_1,\ldots,x_b)$
parametrized by $S$.  The map $E \mapsto \mE$ defines an isomorphism
of functors
$$\ParVect_0(X;x_1,\ldots,x_b) \to \ParBun(X;x_1,\ldots,x_b;\SL(r)) $$
mapping families of semistable bundles to families of semistable
bundles.  There are similar statements in the algebraic category.

We warn the reader that a homomorphism $G \to H$ does not in general
map the Weyl alcove for $G$ into the Weyl alcove for $H$.  This makes
directly associating a morphism of functors to any such homomorphism
problematic, and we will avoid doing so.

\subsection{Equivalence with equivariant bundles}

Parabolic principal bundles are equivalent to bundles equivariant for
a finite group, just as in the vector bundle case.

Let $\Gamma$ denote a finite group acting generically freely on a
curve $\ti{X}$, and let $X = \Gamma \backslash \ti{X}$.  Suppose that
the projection
\label{pisec} $\pi: \ \ti{X} \to X$ has ramification points
$x_1,\ldots,x_b$.  We denote the inverse image of $x_i$ in $\ti{X}$ by
$\ti{x}_i$.  The stabilizer of $\ti{x_i}$ under $\Gamma$ is denoted
$\Gamma_{\ti{x_i}}$.  We fix a generator $\gamma_{\ti{x_i}}$ of
$\Gamma_{\ti{x_i}}$, so that its action in a neighborhood of
$\ti{x_i}$ is given by multiplication by a primitive root of unity.

Let $\ti{E} \to \ti{X}$ be a $\Gamma$-equivariant vector bundle.
Define $E$ to be the vector bundle whose sheaf of sections is sheaf of
$\Gamma$-invariant sections of $\ti{E}$.  The parabolic structures are
the filtrations of $E_{x_i}$ induced by order of vanishing at the
ramification points.  The markings $\mu_i$ are the logarithms of the
eigenvalues of the generator of $\Gamma$, acting on
$\ti{E}_{\ti{x_i}}$, for any $\ti{x}_i$ in the fiber over $x_i$.  Let
$\on{Vect}_\Gamma$ denote the functor which assigns to any complex
manifold $S$, the isomorphism classes of $\Gamma$-equivariant bundles
$\ti{E} \to \ti{X}$.  The map $\ti{E} \mapsto E$ defines an
isomorphism of functors, $\on{Vect}_\Gamma(\ti{X}) \to
\ParVect(X;x_1,\ldots,x_b)$ \cite{ms:pb,fu:se,bo:re,bi:pb}.

Let $\Bun_\Gamma(\ti{X};\ti{x};G)$ be the functor which assigns to any
complex manifold $S$, the isomorphism classes of $\Gamma$-equivariant
principal $G$-bundles $\ti{\mE} \to \ti{X}$.  We will sketch a proof
of the following theorem:

\begin{theorem}  \label{isofunct}
There exists an isomorphism of functors
$$\Bun_\Gamma(\ti{X};\ti{x};G) \to \ParBun(X;x;G) $$
mapping families of semistable bundles to semistable bundles, and a
similar isomorphism in the algebraic category.
\end{theorem}

That is, there is a natural bijection between isomorphism classes of
$\Gamma$-equivariant bundles (resp. semistable bundles) on $S \times
\ti{X}$, and isomorphism classes of parabolic bundles
(resp. semistable bundles) on $S \times X$.  Parabolic bundles with
parabolic weights $\mu_i$ at the points $x_i$ are mapped to
$\Gamma$-equivariant bundles with action at $\ti{x}_i$ in the
conjugacy class given by $\mu_i$.

Let $\ti{\mE} \to \ti{X}$ be a $\Gamma$-equivariant principal
$G$-bundle.  We suppose for simplicity there is a single fixed point
$\ti{x} = \ti{x}_j$ with marking $\mu = \mu_j$ and stabilizer
$\Gamma_{\ti{x}_j} = \Gamma$. \label{Usec} Choose a neighborhood
$\ti{U} \to U$ with local coordinate $z$ so that the projection is
given by $z \mapsto z^N$, and the action of $\Gamma$ by $z \mapsto
\exp(2 \pi i /N)z $.  By the equivariant Oka principle of Heinzner and
Kutzchebauch \cite[Section 11]{hk:eo}, after shrinking $\ti{U}$ we may
assume that $\ti{\mE}$ is $\Gamma$-trivial over $\ti{U}$.  That is,
there exists a $\Gamma$-equivariant biholomorphic map $\tau: \
\ti{\mE}|_{\ti{U}} \to \ti{U} \times G$ such that the action of
$\Gamma$ is given by $\gamma(z,g) = (\exp(2\pi i/N) z, \exp(\mu)g) .$
Consider the one parameter subgroup,
$$ \ \C^* \to G, \ \ z \mapsto z^{N\mu/2 \pi i} : = \exp(\ln(z)N \mu/
2\pi i) .$$
Let $\ti{\Sigma}^{-N\mu}$ denote the set of $\Gamma$-invariant
meromorphic sections $s: \ \ti{U} \to \ti{\mE}$ such that $s(z)
z^{-N\mu/2\pi i}$ is regular on $\ti{U}$.  $ \ti{\Sigma}^{-N\mu}$
contains the section given locally by $s_0(z) = z^{N\mu/2\pi i}$.  

We wish to show that there is a parabolic bundle $(\mE,\varphi,\mu)$
isomorphic to $\Gamma \backslash \ti{\mE}$ over $\Gamma \backslash
(\ti{X} \backslash \ti{x})$ such that $\ti{\Sigma}^{-N\mu}$ is the set
of sections of $\mE$ over $U$.  Form a bundle $\ti{\mE}^{-N\mu}$
by patching together $ \ti{\mE}|_{\ti{X} \backslash \{
\ti{x} \}}$ with $\ti{U} \times G$, using the transition map
$z^{-N\mu/2\pi i}$.  The action of $\Gamma$ extends to
$\ti{\mE}^{-N\mu}$ and is trivial near $x$.  Define $\mE = \Gamma
\backslash \ti{\mE}^{-N\mu}$.  Since $\Gamma$ acts trivially in the
fiber at the ramification point, $\mE$ is a principal $G$-bundle.  Let
$\varphi \in \mE_{x_j}/P_j$ denote the parabolic reduction given as
$P_j$ in the trivialization at $x_j$.  We leave it to the reader to
check that the definition of $(\mE,\varphi)$ is independent of the
choices (this depends essentially on the assumption $\alpha_0(\mu) <
1$) and defines an isomorphism of functors.

We construct a one-to-one correspondence between parabolic reductions
of $\ti{\mE}$ and $\mE$, which maps the degree to a multiple of the
parabolic degree.  Let $\ti{\mE} \to \ti{X}$ be a $\Gamma$-equivariant
bundle, and $\mE = \ti{\mE}^{-N \mu}/ \Gamma$.  Any parabolic
reduction $\sigma: \ X \to \mE/P$ induces a $\Gamma$-invariant
parabolic reduction $\ti{\sigma}$ of $\ti{\mE}$ and vice-versa, since
$G/P$ is complete.  Fix a local trivialization $\ti{U}_i \times G$
near $\ti{x}_i$, so that the action of $\Gamma$ is given by
$\exp(\mu_i)$ on the fiber.  The bundle $\mE$ is formed by twisting by
$z^{-N\mu_i/2\pi i}$ near $\ti{x}_i$, and taking the quotient by
$\Gamma$.  Using this local trivialization, the fixed point set of
$\Gamma$ on $\ti{\mE}_{\ti{x}_i}/P$ has components indexed by the
double coset space of the Weyl group $ W_{P_i} \backslash W / W_P$:
\begin{equation} \label{fpset}
 (\ti{\mE}_{\ti{x}_i}/P)^\gamma \cong (G/P)^{\exp(\mu)} = \bigcup_{w \in
 W_{P_i} \backslash W / W_P} L w P
\end{equation}
where $L$ is the standard Levi subgroup of $P$. \label{levisec} 

\begin{lemma}  $\ti{\sigma}(\ti{x}_i) \in LwP $, if and only if the relative position of
$(\varphi_i,\sigma(x_i))$ is $[w]$. 
\end{lemma}

\begin{proof}  Let $O_w \subset G/P$ be the open cell containing $\Ad(w)P$, that is,
$O_w = w B^{-} P$, where $B^{-}$ is the Borel opposite to $B$.  Let
$C_w = BwP$ be the Schubert cell containing $Ad(w)P$.  The set of
elements $g \in G/P$ in relative position $[w]$ is $C_w$.
\label{Csec} Let $R_-(P)$ be the set of weights of $\g/\p$, i.e. roots
of the negative unipotent complementary to $P$.  Let $ f : \ O_w \to
\times \g_\alpha $ be the $T$-equivariant isomorphism of the open cell
$O_w$ with the product of root spaces $ \g_\alpha $ for ${\alpha \in
wR_-(P) }$.  The image $ f(C_w)$ is the product of $\g_\alpha$ for
$\alpha \in w R_-(P) \cap R_+(B)$.  Therefore it suffices to show that
each component $ \ti{\sigma}^{-N\mu}_\alpha$ is regular at $z = 0$ and
$ \ti{\sigma}^{-N\mu}_\alpha(0) = 0 $ unless $ \alpha \in R_+(B)$.
Since $\ti{\sigma}$ is $\gamma$-invariant, $ \ti{\sigma}_\alpha(\gamma
\cdot z) = \exp( 2 \pi i (\mu,\alpha) ) \ti{\sigma}_\alpha(z) .$ Hence
$$ \ti{\sigma}_\alpha(z) = \sum_{ j\ge 0, \ j \equiv_N (N\mu,\alpha)}
c_j z^j .$$
Therefore
$$ \sigma^{-N\mu}_\alpha(z) = \sigma_\alpha(z) z^{- (N \mu,\alpha)} =
\sum_ { j \ge 0, \ j \equiv_N (N \mu,\alpha) } c_j z^{j -
(N\mu,\alpha)} .$$
It follows that $\sigma^{-N\mu}_\alpha(z) = 0 $ at $z = 0$ if
$(\alpha,\mu) < 0$, and is regular in any case.
\end{proof}

We compute the parabolic degree of $\sigma$ as follows.  In the local
trivialization of $\ti{\mE}$ near $\ti{x}_i$, the reduction
$\ti{\sigma}$ is given by $\ti{\sigma}_i(z) P$ for some map
$\ti{\sigma}_i: \ \ti{U}_i \to G$.  By the previous paragraph we may
assume $\ti{\sigma}_i(0) = n_i$, for some representative $n_i$ of
$w_i$ the relative position of $\sigma(x_i)$ and $\varphi_i$.  By
equivariant Oka \cite{hk:eo} applied to the $\Gamma$-equivariant
$P$-bundle corresponding to $\sigma$, we may assume that
$\ti{\sigma}_i(z) = n_i$ is constant.  The bundle
$(\ti{\sigma}^{-N\mu})^* \ti{\mE}^{-N\mu}$ is formed by gluing
$\ti{\sigma}^* \ti{\mE} \backslash \bigcup \ti{\sigma}^*
\ti{\mE}_{x_i}$ with $\bigcup \ti{\sigma}^* \ti{\mE} |_{\ti{U}_i}$
using the maps $ \Ad(n_i) z^{-N\mu_i / 2\pi i} = z^{-N w_i \mu_i/2\pi
i}$.  This implies that the gluing maps for $(\ti{\sigma}^{-N\mu})^*
\ti{\mE}^{-N\mu}(\omega_P)$ are $\chi_P( \Ad(n_i) z^{-N\mu_i/2\pi i})
= z^{- \omega_P(N w_i \mu_i)}.$ The degree of the line bundle is
therefore
$$ \deg((\ti{\sigma}^{-N\mu})^* \ti{\mE}^{-N \mu}(\omega_P)) = \deg
(\ti{\sigma}^* \ti{\mE}(\omega_P)) - \sum_{i = 1}^b N \omega_P(w_i
\mu_i) .$$
Since $\sigma = \Gamma \backslash \ti{\sigma}^{N \mu}$, the degree of
$ \sigma^* \mE(\omega_P)$ is $1/N$ times the degree of $(\ti{\sigma}^{-N
\mu})^* \ti{\mE}^{-N \mu}(\omega_P)$.  Hence,
$$ \deg(\sigma^*\mE(\omega_P)) + \sum_{i = 1}^b \omega_P(w_i \mu_i) =
\frac{1}{N} \deg(\ti{\sigma}^* \ti{\mE}(\omega_P)) .$$
That is, the parabolic degree of $\sigma$ is $ \deg(\ti{\sigma})/N .$

\subsection{Modifications for the algebraic case}

To prove the correspondence in the algebraic category one has to
replace the equivariant Oka principle by a non-abelian cohomology
argument, and the gluing by formal gluing.

\begin{lemma} Let $\ti{\mE} \to \ti{X} \times S$ be a
$\Gamma$-equivariant principal $G$-bundle.  For any $s \in S$, there
exists a neighborhood that is the product of an \'etale neighborhood
in $S$ and formal neighborhood of $\ti{x}$ in $\ti{X}$, such that the
action of $\Gamma$ on the restriction of $\ti{\mE}$ is of product
form.
\end{lemma}

\begin{proof} 
Over the formal disk $D = \Spec(\C[[z]])$ at $\ti{x}$, the bundle
$\ti{\mE}$ is trivial and the action of $\Gamma$ is given by $
\gamma(z,\zeta) = (\gamma z, g(\gamma,z) \zeta) $ for some $g: \Gamma
\to G[[z]]$.  Since $\gamma^N =1 $ we have
$$ g(\gamma,\gamma^{N-1} z) g(\gamma,\gamma^{N-2} z) \ldots g(\gamma,z) =
e .$$ 
In particular, $g(\gamma,0)^N = e$.  More generally, if $\ti{\mE}$ is a
$\Gamma$-equivariant bundle over $D_R := \Spec(R[[z]])$, where $R$ is
any $\C$-algebra, then the action is given by an automorphism $g \in
G(R[[z]])$.  Let $\ti{\mE}$ be a $\Gamma$-equivariant bundle over $D_R$.
We wish to show that there exists an automorphism $\tau \in G(R[[z]])$
which transforms the $\Gamma$-action on $\ti{\mE} |_{D_R} \cong D_R
\times G$ to the product action, that is, $ \tau(\gamma z)^{-1}
g(\gamma,z) \tau(z) = g(\gamma,0) .$ Consider the element of
$C^1(\Gamma,G(R[[z]]))$ defined by $\gamma \mapsto g(\gamma,z)$.
Since
$$g(\gamma_1\gamma_2,z) = g(\gamma_1,\gamma_2z) g(\gamma_2z) =
\gamma_2^* g(\gamma_1,z) g(\gamma_2,z) $$
$g(\cdot,z)$ is a cocycle in the cohomology of $\Gamma$ with values in
$G[[z]]$. Similarly $g(\cdot,0) \in Z^1(\Gamma,G(R))$ which maps to
$Z^1(\Gamma,G(R[[z]]))$.  We claim there exists a $0$-chain $\tau$
such that $(\delta \tau): \ g(\cdot,z) \mapsto g(\cdot,0)$.  We
construct $\tau$ order-by-order.  Let $G_l = G(R[z]/z^{l}= 0)$.  Let
$N_l$ be the kernel of the truncation map $G_{l+1} \to G_l$.  The
exact sequence of groups
$$ 1 \to N_l \to G_{l+1} \to G_l \to 1 $$
induces an exact sequence of pointed sets in non-abelian cohomology
(see e.g. \cite[p. 49]{se:ga})
$$ H^1(\Gamma,N_l) \to H^1(\Gamma,G_{l+1}) \to H^1(\Gamma,G_l) .$$
Since $N_l$ is a nilpotent, $H^1(\Gamma,N_l)$ is trivial, by induction
on the length of the central series which reduces to the case that
$N_l$ is a $\Gamma$-module.  Therefore, $H^1(\Gamma,G_{l+1})$ injects
into $H^1(\Gamma,G_l)$ for all $l$.  The complexes
$C^{0}(\Gamma,G_l),C^{1}(\Gamma,G_l)$ satisfy the Mittag-Leffler
condition: the image of $C^0(\Gamma,G_{l'})$
(resp. $C^{1}(\Gamma,G_{l'})$) in $C^{0}(\Gamma,G_l)$
resp. $C^{1}(\Gamma,G_l)$ stabilizes as $l' \to \infty$.  Indeed, let
$f_l : \ \Spec(R_l) \to G$.  Extending $f_l$ to a map $f_{l+1}: \
\Spec(R_{l+1}) \to G$ is equivalent to extending the map $f_0^{-1}
f_{l+1} $; the latter extends because $G$ is isomorphic to $\g$ near
the identity.  Therefore, $G_{l+1} \to G_l$ is surjective, which
implies the same result for the chain complexes.  The Mittag-Leffler
condition implies that
$$ H^1(\Gamma,G[[z]]) = \lim_{l \to \infty} H^1(\Gamma,G_l) $$ 
(see \cite[II.9.1]{ha:al} for the abelian case) and therefore also
injects into $H^1(\Gamma,G)$.  The claim follows since $g(\cdot,z)$
and $g(\cdot,0)$ both map to $g(\cdot,0)$ in $H^1(\Gamma,G(R))$. 
\end{proof}

Recall the description of bundles on $\ti{X}$ by formal gluing data
\cite{be:de}, \cite[Section 3]{la:li}.  For any algebra $R$, let
$\ti{X}_R := \ti{X} \times \Spec(R) .$ Let $T$ denote the functor from
algebras to sets which associates to $R$ the set of isomorphism
classes of triples $(\ti{\mE},\rho,\sigma)$, where $\ti{\mE}$ is a
$G$-bundle over $\ti{X}_R$, $\rho$ is a trivialization over $(\ti{X}
\backslash \{ x \})_R$, and $\sigma$ is a trivialization over the
formal disk $D_R$.  Then $T$ is represented by $G(R((z)))$
\cite[3.8]{la:li}.  Choose a set of trivializations of $\ti{\mE}$ in
formal neighborhoods of the form $D_R$ as described above.  Let
$\ti{\mE}^{- N\mu}_R$ denote the bundle obtained from twisting by $z^{-N
\mu/2 \pi i} \in G(R((z)))$.  The bundles $\ti{\mE}^{-N \mu}_R$ are
canonically isomorphic away from $\ti{x}$, the canonical isomorphisms
extend to $\ti{X}$, and the extension preserves the parabolic
structures at the ramification points.  By a simple case of \'etale
descent, the bundles $\ti{\mE}^{-N\mu}_R$ patch together to a bundle
$\ti{\mE}^{-N\mu} \to S \times \ti{X}$.  Since the gluing data for
$\ti{\mE}^{-N \mu}$ are $\Gamma$-invariant, they define a $G$-bundle $\mE
\to S \times X$ with parabolic structure.

\subsection{Canonical reductions}

If a parabolic vector bundle $E$ is unstable, the Harder-Narasimhan is
a canonical sequence of sub-bundles violating the semistability
condition.  There is a unique sub-bundle $E_{1} \subset E$
\label{maxsec} such that the slope $\mu(E_{1})$ is maximal among
all sub-bundles, and the rank of $E_{1}$ is maximal among sub-bundles
with that slope.  The Harder-Narasimhan filtration 
$$E_\bullet = \{ E_1 \subset \ldots \subset E_k = E \}$$ 
is defined inductively by $E_{i+1}/E_i = (E/E_i)_{1} .$ It follows
from the definition that the quotients $E_{i+1}/E_i$ of the canonical
filtration are semistable, the slopes $\mu_i = \mu(E_{i}/E_{i-1})$
decreasing, and $E_{\bullet}$ is the unique filtration with slopes
$\mu_i$ and ranks $r_i = \dim(E_{i})$.

Atiyah-Bott \cite[Section 10]{at:mo} construct a canonical parabolic
reduction $\sigma_E: \ X \to \mE/P$ generalizing the Harder-Narasimhan
filtration, using the adjoint bundle $\mE(\g)$ of $\mE$.  \label{canred}
The Harder-Narasimhan filtration 
$$\ldots \mE(\g)_{-1} \subset \mE(\g)_0 \subset \ldots $$ 
of $\mE(\g)$ has a term $\mE(\g)_0$ such that $\mE(\g)_0 /
\mE(\g)_{-1}$ has degree zero, and $\mE(\g)_0$ is a bundle of
parabolic Lie algebras for some parabolic subgroup $P$ and defines a
reduction $\sigma_{\mE}: \ X \to \mE/P$. The value $\sigma_{\mE}(x)$
of the reduction at $x \in X$ is the unique fixed point for
$\mE(\g)_0$ acting infinitesimally on the fiber $\mE_x/P$.

The canonical reduction $\sigma_{\mE}$ is functorial for homomorphisms
$\phi: \ G \to G'$ such that the associated Lie algebra map $D\phi: \
\g \to \g'$ is injective.  That is, for any principal $G$-bundle $\mE
\to X$, there is a parabolic subgroup $P'$ of $G'$ such that $P$ is
the inverse image of $P'$ under $\phi$ and $\sigma_{\mE}$ is the
inverse image of $\sigma_{\phi_*E}$ under the map $\mE/P \to \mE'/P'$.
Indeed, the image of $\mE(\g)_0$ in $\phi_* \mE(\g')$ is contained in
$\phi_*\mE(\g')_0$, for reasons of degree, and maximality of
$\mE(\g)_0$ among sub-bundles with the same degree implies that
$\mE(\g)_0$ contains the inverse image of $\phi_*\mE(\g')_0$.  Hence
$\p = D\phi^{-1}(\p')$.

Define a notion of slope for parabolic reductions as follows.  Let
$\Lambda_P^*$ denote
\label{charsec}
the abelian group of weights of characters of $P$, and $\Lambda_P$ its
dual.  For a principal $G$-bundle $\mE \to X$ and parabolic reduction
$\sigma: \ X \to \mE/P$, the {\em slope} $\mu(\sigma) \in \Lambda_P$ is
given by
$$\mu(\sigma): \lambda \mapsto \deg(\sigma^*\mE(\lambda)) $$  
for $\lambda \in \Lambda_P^*$.  \label{assocsec} The {\em type} of $\mE$
is the slope $\mu(\sigma_{\mE})$ of its canonical reduction.
$\mu(\sigma_{\mE})$ lies in the interior of the open face of $\t_+$
corresponding to $P$.

\begin{lemma}
The canonical reduction $\sigma_{\mE}$ is the unique parabolic reduction
with slope $\mu(\sigma_{\mE})$.
\end{lemma}

\begin{proof}  Consider an embedding $\phi: \ G \to
\Gl(V)$, and let $\sigma_{\phi_*E}: \ X \to \phi_* \mE/ P'$ denote its
canonical reduction.  Let $\sigma: \ X \to \mE/P$ be another reduction
with slope $\mu(\sigma_{\mE})$, and $\phi_* \sigma$ be the parabolic
reduction of $\phi_* \mE$ to $P'$ induced by $\sigma$.  Since
$\deg((\phi_* \sigma)^* \phi_* \mE(\lambda)) =
\deg(\sigma^*\mE(D\phi^*\lambda))$ for any weight $\lambda \in
\Lambda_{P'}^*$, $\mu(\sigma_{\phi_*E}) = \mu(\phi_* \sigma) .$ Since
the Harder-Narasimhan filtration is the unique filtration of its
slope, $\phi_* \sigma = \phi_* \sigma_{\mE}$.  This implies that $\sigma =
\sigma_{\mE}$.
\end{proof}

Using the equivalence of parabolic bundles with equivariant bundles
one can extend the theory of the canonical reduction to parabolic
bundles.  Let $\gamma: \ X \to X$ be an automorphism of the curve $X$.
If $\mE$ is an $\gamma$-equivariant bundle, then the canonical reduction
is $\gamma$-invariant, since it is the unique reduction with its
slope.  Let $\Gamma$ be a group of automorphisms of $\mE$. We will call
{\em $\Gamma$-stable} (resp. {\em $\Gamma$-semistable}) if
\begin{equation} \label{ss}
 \deg (\sigma^*\mE(\lambda)) \leq 0 \text{ (resp. } < 0 \text{ ) }
\end{equation}
for all $\Gamma$-invariant parabolic reductions $\sigma: \ X \to \mE/P$
and weights $\lambda \in \Lambda_P^*$.  Since the canonical reduction
is the unique reduction of its slope, a principal $G$-bundle is
$\Gamma$-semistable if and only if it is ordinary semistable.  On the
other hand, $\Gamma$-stability, or $\Gamma$-irreducibility of $\mE$ is
not in general the same as ordinary stability or irreducibility.  For
any parabolic bundle $\mE = (\mE,\{ (\varphi_i,\mu_i) \})$, let $\sigma_{\mE}$
denote its canonical reduction, defined by the one-to-one
correspondence between invariant parabolic reductions of $\ti{\mE}$ and
parabolic reductions of $\mE$.  Define the slope of a parabolic
reduction $\sigma: \ X \to \mE/P$ by
$$ \mu(\sigma): \ \lambda \mapsto \deg \sigma^* \mE(\lambda) + \sum
\lambda( w_i \mu_i) .$$
The {\em type} of $\mE$ is the slope of $\sigma_{\mE}$; by the discussion
above $\sigma_{\mE}$ is the unique reduction of this slope.

\subsection{Grade equivalence}

The rest of this section is included for the sake of completeness,
and is not needed for the main result. 

We extend Ramanathan's notion of grade equivalence to parabolic
bundles.  First, let $\mE \to X$ be a $G$-bundle, and $\sigma: X \to
E/P$ be a parabolic reduction.  Let \label{rsec} $r: \ P \to L$ the
projection to a Levi subgroup $L \subset P$, and $\iota: \ L \to G$
the inclusion of $L$ in $G$.  The reduction $\sigma$ is {\em
admissible} if $\deg(\sigma^*\mE(\lambda)) = 0$ for all weights
$\lambda$.  The equivalence relation on semistable bundles generated
by
$$\mE \sim \iota_* r_* \sigma^* \mE,$$ 
as $\sigma$ ranges over all admissible reductions, is called {\em
grade equivalence} \cite{ra:th,ra:th2}.  For any semistable bundle
$\mE \to X$, there is a semistable bundle $\Gr(\mE)$, unique up to
isomorphism, defined by the condition that there is an admissible
reduction $\sigma: \ X \to \mE/P$ such that $r_* \sigma^* \mE $ is
stable and $\Gr(\mE) \cong \iota_*r_*\sigma^*E$.  The set of
isomorphism classes of semistable $G$-bundles $\mE$ such that $\mE
\cong \Gr(\mE)$ form a set of representatives for the equivalence
classes of semistable $G$-bundles over $X$.  That is, two bundles
$\mE_1,\mE_2 \to X$ are grade equivalent, if and only if their grade
bundles $\Gr(\mE_1),\Gr(\mE_2)$ are isomorphic.

For equivariant bundles we define grade equivalence to be the
equivalence relation generated by $\mE \sim \iota_* r_* \sigma^* \mE$,
where $\sigma: \ X \to \mE/P$ is $\Gamma$-invariant.  

Let $ \mE \to S \times X $ be a family of $\Gamma$-equivariant principal
$G$-bundles, and $\mE_0$ a bundle such that $\mE_s$ is $\Gamma$-isomorphic
to $\mE_0$ for $s$ varying in a dense open subset of $S$.  The
equivalence relation generated by $\mE_s \to \mE_0$ for any $s \in S$ is
called $S$-equivalence.  By \cite[Proposition 3.24]{ra:th}, this is
the same as grade-equivalence.

To define grade equivalence for parabolic bundles, first let $E$ be a
parabolic vector bundle, with filtration $E_{x_i}^\bullet$.  If $F
\subset E$ is any sub-bundle, the filtrations $E_{x_i}^j$ induce
filtrations on the fibers of the graded bundle $F \oplus E/F$ at the
points $x_i$, and we say that $F \oplus E/F$ is parabolic grade
equivalent to $E$.  

This construction generalizes to arbitrary type as follows.  let
$(\mE,\varphi_1, \ldots,\varphi_b,\mu_1,\ldots,\mu_b)$ be a parabolic
bundle, $\sigma: \, X \to \mE/P$ a parabolic reduction, and $r_*
\sigma^* \mE$ the associated $L$-bundle.  Let $\Aut(\mE_{x_i})$ denote
the group of $G$-equivariant automorphisms of $\mE_{x_i}$.
$\Aut(\mE_{x_i})$ is isomorphic to $G$, and the stabilizer $P' =
\Aut(\mE_{x_i})_{\sigma(x_i)}$ is isomorphic to $P$.  Similarly, the
stabilizer $P_i' = \Aut(\mE_{x_i})_{\varphi_i}$ is isomorphic to
$P_i$.  In the vector bundle case, $P_i'$ is the group of
automorphisms preserving the filtration $E_{x_i}^\bullet$.  The
intersection $P' \cap P_i'$ is a subgroup isomorphic to $w_i P_i \cap
P$, where $w_i$ is the relative position of $\varphi_i$ and
$\sigma(x_i)$.  Its image in $\Aut(r_* \sigma^* \mE_{x_i})$ is a
parabolic subgroup, isomorphic to $w_i P_i \cap L$.  Therefore, it has
a unique closed orbit in $r_* \sigma^* \mE_{x_i}/ (w_i P_i \cap L)$.
Since $r_* \sigma^* \mE_{x_i}/ (w_i P_i \cap L)$ injects into $\iota_*
r_* \sigma^* \mE_{x_i}/P_i$, we get a reduction of $\iota_* r_*
\sigma^* \mE_{x_i}/P_i$ we denote by $\iota r \sigma \varphi_i$.  Let
{\em parabolic grade equivalence} be the equivalence relation
generated by
$$(\mE,\varphi_1, \ldots,\varphi_b,\mu_1,\ldots,\mu_b) \sim (\iota_* r_*
\sigma^*\mE,\iota r \sigma \varphi_1, \ldots, \iota r \sigma \varphi_b,
\mu_1,\ldots,\mu_b ) .$$
We claim this equivalence relation corresponds to grade equivalence
for equivariant bundles.  Let $\ti{\mE} \to \ti{X}$ be a
$\Gamma$-equivariant bundle, and $(\mE,\varphi_1,\ldots,\varphi_b,
\mu_1,\ldots,\mu_b )$ the corresponding parabolic bundle.  Let
$\ti{\sigma}$ be a $\Gamma$-invariant parabolic reduction to a
parabolic subgroup $P$ and $\sigma$ the corresponding parabolic
reduction of $\mE$.  Let $\ti{U}_i \times G$ be a local trivialization
near $x_i$, so that the action of $\Gamma$ is $(z,g) \mapsto (\exp
(2\pi i /N)z, \exp(\mu_i)g)$ and $\ti{\sigma}(z) = w_i^{-1} P$, for
some $w_i \in W_P \backslash W$.  The local trivialization of
$\ti{\sigma}^*\ti{\mE}$ induces a local trivialization of $r_*
\ti{\sigma}^* \ti{\mE}$ near $x_i$.  The action of $\Gamma$ is given in
this local trivialization by
$$ (z, l) \mapsto (\exp (2\pi i /N)z, \exp(w_i \mu_i) l ) .$$
Let $\mu_{L,i}$ be the unique point in the positive chamber for $L$
conjugate to $w_i \mu_i$.  The parabolic bundle corresponding to $r_*
\ti{\sigma}^* \ti{\mE}$ is $(\mE_L,\varphi_L,\mu_L)$ where $\mE_L = (r_*
\ti{\sigma}^* \ti{\mE})^{-N \mu_L}/\Gamma$.  Define $w_{L,i} \in W_L$
the Weyl group for $L$ by $ w_{L,i} \mu_{L,i} = w_i \mu_i .$ By
definition the gluing maps for $ (r_* \ti{\sigma}^* \ti{\mE})^{-N
\mu_L}$ are given $w_{L,i} z^{N\mu_{L,i}/2 \pi i}$.  The gluing map
for $r_* (\ti{\sigma}^{-N\mu})^* \ti{\mE}^{-N\mu} $ is $ z^{-N w_i
\mu_i/ 2\pi i}$.  Since $ z^{-Nw_i \mu_i/2 \pi i} w_{L,i}
z^{N\mu_{L,i}/2 \pi i } = w_{L,i} $ is regular at $z =0$, the bundles
$(r_* \ti{\sigma}^* \ti{\mE})^{-N \mu_L}$ and $ r_*
(\ti{\sigma}^{-N\mu})^* \ti{\mE}^{-N\mu} $ are isomorphic.  Therefore,
their quotients by $\Gamma$ are isomorphic.  The parabolic structure
for $(r_* \ti{\sigma}^* \ti{\mE})^{-N \mu_L}$ at $x_i$ is $r(P \cap w_i^{-1} P_i)$ in the trivialization near
$x_i$. This completes the proof of the claim.

\subsection{Coarse moduli spaces}

Let $\ol{\Bun}^{\ss}(X)$ denote the functor which associates to any
scheme $S$ the set of grade equivalence classes of semistable
algebraic principal $G$-bundles over $S \times X$.  The main result of
Ramanathan's thesis \cite{ra:th2} (see also \cite{fa:st})) is the
existence of an irreducible, normal projective variety $\M_G(X)$ and a
morphism $\ol{\Bun}^{\ss}(X) \to \Hom(\cdot, \M_G(X))$ that is a
coarse moduli space for $\ol{\Bun}^{\ss}(X)$. \label{MGsec} By
definition, a {\em coarse moduli space} for a functor $F$ is a scheme
$M$ and a morphism $\rho: \ F \to \Hom( \cdot, M)$ such that (i)
$\rho$ induces a bijection of points $\rho(*): \ F(*) \to \Hom(*,M)$,
where $* = \Spec(\C)$, and (ii) for any scheme $N$ and morphism $\chi:
\ F \to \Hom(\cdot, N)$, there is a unique morphism $\phi: \
\Hom(\cdot, M) \to \Hom(\cdot,N) $ such that $\chi = \phi \circ \rho$.
Usually, we omit the morphism $\rho$ from the notation.  

Let $\ol{\Bun}^{\ss}_\Gamma(\ti{X})$ denote the functor that assigns
to any scheme (or complex manifold, in the analytic category) $S$ the
set of grade-equivalence classes of $\Gamma$-equivariant bundles over
$S \times \ti{X}$.

\begin{theorem} \label{coars1}
There is a normal projective variety $\M_{G,\Gamma}(\ti{X})$ that is a
coarse moduli space for $\ol{\Bun}^{\ss}_\Gamma(\ti{X})$.
\end{theorem}
\noindent Sketch of Proof: We realize $\M_{G,\Gamma}(\ti{X})$ as a
subquotient of the moduli space of bundles with level structure.
Recall that a {\em level structure} on $\mE$ at a point $y \in \ti{X}$
is a point $e_y$ in the fiber $\mE_y$.  Bundles with level structure
have no automorphisms, since the map $\Aut(\mE) \to \Aut(\mE_y)$ is
injective.  A {morphism of bundles with level structure}
$(\mE_1,e_{1,y}),(\mE_2,e_{2,y})$ is a morphism $\varphi: \ \mE_1 \to
\mE_2$ such that $\varphi(e_{1,y}) = e_{2,y}$.  Let $\Bun(\ti{X};y_1,
\ldots,y_m;G)$ denote the functor which associates to any scheme $S$
the set of isomorphism classes of $G$-bundles over $S \times \ti{X}$
with level structures at points $y_1,\ldots,y_m \in \ti{X}$.  Let
$\Bun^{\on{ss}}(\ti{X};y_1, \ldots,y_m;G)$ denote the open subfunctor
defined by the condition that the underlying bundle is semistable.
$\Bun^{\on{ss}}(\ti{X},y_1,\ldots,y_m;G)$ is represented by a smooth
quasi-projective \label{MGYsec} moduli space
$\M_G(\ti{X},y_1,\ldots,y_m)$, see \cite[Part 4]{se:fi},\cite{hu:pa};
for arbitrary $G$ one needs the embedding arguments in
\cite[4.8.1]{ra:th}.  The right action of $G$ on the fiber at each
marked point induces an action of $G^m$ on
$\M_G(\ti{X},y_1,\ldots,y_m)$, with good quotient $\M_G(\ti{X})$.
Using Hilbert schemes as in \cite[Section 5]{ra:th2} one may construct
a universal space $\U_G(\ti{X},y_1,\ldots,y_m)$ for $G$-bundles with
level structure at $y_1,\ldots,y_m$, such that
$\U_G(\ti{X},y_1,\ldots,y_m) \to \M_G(\ti{X},y_1,\ldots,y_m)$ is a
good quotient.

Suppose that the set $\{ y_1,\ldots,y_m \}$ is invariant under
$\Gamma$, and the stabilizers $\Gamma_{y_i}$ are trivial.  An {\em
equivariant bundle with level structure} is an equivariant bundle
$\mE$ with level structure $e_{y_1},\ldots,e_{y_m}$ such that
$\gamma(e_{y_i}) = e_{\gamma(y_i)}$.  Let
$\M_{G,\Gamma}(\ti{X},y_1,\ldots,y_m)$ denote the set of isomorphism
classes of equivariant bundles with level structure whose underlying
bundle is semistable.  Since bundles with level structure have no
automorphisms, forgetting the equivariant structure defines an
injection
$$\M_{G,\Gamma}(\ti{X},y_1,\ldots,y_m) \to \M_G(\ti{X},y_1,\ldots,y_m)
.$$
The image is the fixed point set of the action of $\Gamma$, which is a
smooth quasi-projective variety.  Let $G_\Gamma^m$ denote the subgroup
of $G^m$ invariant under the action of $\Gamma$ on $G^m$ induced by
the action of $\Gamma$ on the set $\{y_1,\ldots,y_m \}$.  An
observation of Ramanathan is that if $f: X \to Y$ is an affine
morphism of $G$-varieties, and $Y$ has a good quotient, then so does
$X$ \cite[3.12]{ne:mo}.  Note that
$$ \M_{G,\Gamma}(\ti{X}) \times_{G_\Gamma^m} G^m \to \M_G(\ti{X})
\times_{G_\Gamma^m} G^m \to \M_G(\ti{X}) $$
are affine $G$-morphisms; it follows that the action of $G_\Gamma^m$
on $\M_{G,\Gamma}(\ti{X},y_1,\ldots,y_m)$ has a good quotient, which
we denote $\M_{G,\Gamma}(\ti{X})$.  A good quotient is a categorical
quotient, hence $\M_{G,\Gamma}(\ti{X})$ is normal.

We will show that $\M_{G,\Gamma}(\ti{X})$ is a coarse moduli space for
the functor of equivalence classes of $\Gamma$-equivariant bundles.
Let $\mE$ be a $\Gamma$-equivariant semistable bundle over $S \times
\ti{X}$, and $s$ any point in $S$.  In a neighborhood $S_1$ of $s$,
$\mE$ admits equivariant level structures at $y_1,\ldots,y_m$ and
defines $S_1 \to \M_{G,\Gamma}(\ti{X},y_1,\ldots,y_m)$. If $\mE_s$ are
equivariantly isomorphic for $s$ in an open subset $S_0 \subset S$,
then the image of $S_0 \cap S_1$ in
$\M_{G,\Gamma}(\ti{X},y_1,\ldots,y_m)$ is contained in the closure of
a single orbit.  Conversely, if $[\mE_0] \in
\M_{G,\Gamma}(\ti{X},y_1,\ldots,y_m)$ lies in the closure of the orbit
of $[\mE_1] \in \M_{G,\Gamma}(\ti{X},y_1,\ldots,y_m)$, then forgetting
the level structure shows that $\mE_0$ and $\mE_1$ are equivalent.
Hence the points of $\M_{G,\Gamma}(\ti{X})$ are equivalence classes of
semistable bundles.  For any family $\mE \to S \times \ti{X}$ of
equivariant semistable $G$-bundles which admits equivariant level
structure over $S_1 \subset S$, let $\varphi_{\mE,S_1}: \ S_1 \to
\M_{G,\Gamma}$ denote the map induced by adding some level structure,
$S_1 \to \M_{G,\Gamma}(\ti{X},y_1,\ldots,y_m)$, and then composing
with the projection.  It is clear that $\varphi_{\mE}$ does not depend
on the choice of level structure, so that $\varphi_{\mE,S_1}$ patches
together to a map $\varphi_{\mE}$, and $\mE \mapsto \varphi_{\mE}$
defines a morphism of functors
$$\rho_\Gamma: \ \ol{\Bun}^{\ss}_\Gamma(\ti{X};G) \to \M_{G,\Gamma}(\ti{X}) .$$
Part (ii) of the definition of the definition coarse moduli space
follows from the properties of $\U_{G}(\ti{X},y_1,\ldots,y_m)$ as in
\cite[4.5]{ra:th2}.

Let $\L_G(\ti{X},V) \to \M_{G}(\ti{X})$ be the determinant line bundle
associated to a faithful representation $V$ of $G$, see \cite{be:pic}.
This is an ample line bundle; let $\L_{G,\Gamma}(\ti{X},V)$ denote its
pull-back under the forgetful morphism
$$ f: \M_{G,\Gamma}(\ti{X}) \to \M_G(\ti{X}) .$$
We claim that $\L_{G,\Gamma}(\ti{X},V)$ is ample.  Indeed
$\Hom(\Gamma,L)/L$ is finite for finite $\Gamma$ and linear algebraic
$L$; this is essentially a result of A. Weil \cite{we:coh}, see
Slodowy \cite{sl:two}.  By Zariski's main theorem \cite[4.4]{ega3},
any proper morphism with finite fibers is a finite morphism.
\footnote{Mumford \cite[p.124]{mu:red} credits this result to
Chevalley.}  By \cite[6.6]{ega2}, the pull-back of an ample line
bundle under a finite morphism is ample.  This completes the proof of
the claim.  We remark that in the case $\M_{G,\Gamma}(\ti{X})$ is
smooth, the claim follows from Kodaira's theorem.  By the
correspondence theorem in Section 4, $\M_{G,\Gamma}(\ti{X})$ is
compact.  It follows that $\M_{G,\Gamma}(\ti{X})$ is projective.

Let $\M_G(X;x;\mu) := \M_G(X;x_1,\ldots,x_b,\mu_1,\ldots,\mu_b)$ be
the moduli space \label{modsec} of equivalence classes of parabolic
$G$-bundles on $(X;x_1,\ldots,x_b)$ with markings
$\mu_1,\ldots,\mu_b$.  By the equivalence with equivariant bundles
this is a normal projective variety and a coarse moduli space for the
functor $\ParBun(X;x;\mu;G)$ of grade-equivalence classes of
semistable parabolic bundles with markings $\mu_1,\ldots,\mu_b$.

\section{Narasimhan-Seshadri correspondence for parabolic
$G$-bundles}

In this section we prove the correspondence between flat $K$-bundles
and semistable holomorphic $G$-bundles, for markings $\mu_i$
satisfying $\alpha_0(\mu_i) < 1$.  A related result for projective
varieties $X$ of any dimension is proved in \cite{ba:pr}.  A different
approach to this correspondence in the case $SU(r)$ has been given by
Simpson \cite{si:ha}.

The moduli space of flat $K$-bundles on a punctured surface can be
constructed as in Atiyah-Bott as a symplectic quotient of the affine
space of connections by the gauge group.  Let $X$ be a compact
oriented surface with boundary $\partial X$. Since $K$ is
$1$-connected, any principal $K$-bundle on $X$ is trivial.  Let
$$\A(X) := \Omega^1(X,\k), \ \ K(X) := \Map(X,K) $$
be the space of connections on $X \times K \to X$ and gauge group for
$X \times K$.  Choose an invariant inner product $\Tr( \ , \ ): \ \k
\times \k \to \R$ on $\k$.  The affine space $\A(X)$ has a symplectic
form
$$a_1,a_2 \mapsto \int_X \Tr(a_1 \wedge a_2) .$$
The action of $K(X)$ on $\A(X)$ is Hamiltonian, with moment map given
by the curvature plus restriction to the boundary
$$\A(X) \mapsto \Omega^2(X,\k) \oplus \Omega^1(\partial X,\k), \ \ \ A
\mapsto (F_A,\iota^*_{\partial X} A) .$$
The symplectic quotients of $\A(X)$ by $K(X)$ may be identified with
moduli space of flat connections on $\A(X)$, with fixed holonomy
around the boundary \cite{me:lo}.  Let $b$ denote the number of
components of $\partial X$.  The orbits of $K(X)$ on
$\Omega^1(\partial X,\k)$ are parametrized by $b$-tuples $\mu =
(\mu_1,\ldots,\mu_b) \in \Alc^b$.  Let
$$\Hol_i: \ \Omega^1(\partial X, \k) \to K$$
denote the holonomy around the $i$-th boundary component.  Then two
connections $A_1,A_2 \in \Omega^1(\partial X, \k)$ are in the same
orbit of $K(\partial X)$ if and only if $\Hol_i(A_1)$ is conjugate to
$\Hol_i(A_2)$, for $i = 1,\ldots, b$.  The symplectic quotient
$$ \RR_K(X;\mu_1,\ldots,\mu_b) = \A(X) \qu_\mu \ K(X) = \{ A \in
\A(X), \ F_A = 0, \ \ \Hol_i(A) \in \cC_{\mu_i} \} / K(X) $$
\label{repsec} 
is the moduli space of flat connections on $X \times K$, with fixed
holonomy.  Up to symplectomorphism $\RR_K(X;\mu) :=
\RR_K(X;\mu_1,\ldots,\mu_b)$ does not depend on the choice of marked
points $x_i$, which justifies dropping them from the notation.

The spaces $\RR_K(X;\mu)$ may be identified with moduli spaces of
representations of $\pi_1(X)$ in $K$.  Any flat connection $A$
determines a holonomy representation $\Hol(A): \ \pi_1(X) \to K .$ The
$i$-th boundary component $(\partial X)_i$ determines a conjugacy
class $[(\partial X)_i] \subset \pi_1(X)$.  Two flat bundles are
isomorphic if and only if their holonomy representations are conjugate
by the action of $K$.  Therefore, there is a bijection
\begin{equation} \label{holdes}
 \RR_K(X;\mu) \to \{ \rho \in \Hom(\pi_1(X),K), \ \
\rho([(\partial X)_i]) \subset \cC_{\mu_i} \}/K.\end{equation}

Now suppose that $X$ is a compact, oriented two-manifold without
boundary, and $x_1,\ldots,x_b \in X$ distinct marked points.  Then the
moduli space of flat bundles on $X \backslash \{ x_1, \ldots, x_b \}$
with holonomy around $x_i$ in $\cC_i$ is
$\RR_K(X',\mu_1,\ldots,\mu_b)$, where $X'$ is the manifold obtained by
removing a small open disk containing each marked point $x_i$.  We
denote this space by $\RR_K(X;\mu)$.  In the case $X = \P^1$, the
fundamental group of $ X \backslash \{ x_1,\ldots,x_b \}$ is
$$\pi_1 (X \backslash \{ x_1,\ldots,x_b \}) = <
c_1,\ldots,c_b > / \Pi_{i=1}^b c_i = 1 .$$
By \eqref{holdes}, the moduli space of flat bundles is given by
$$ \RR_K(\P^1;\mu_1,\ldots,\mu_b) = \{ (k_1,\ldots,k_b) \in \cC_{\mu_1}
\times \ldots \times \cC_{\mu_b} \ | \ \Pi_{i=1}^b k_i = e \ \} / K .$$

The moduli spaces on the punctured surface are homeomorphic to moduli
spaces of $\Gamma$-invariant flat connections on a ramified cover
$\ti{X}$.  Suppose $\Gamma$ acts on $ \ti{\mE} := \ti{X} \times K$, so
that the generator of $\Gamma_{\ti{x}_i}$ acts on the fiber
$\ti{\mE}_{\ti{x}_i}$ by an element in the conjugacy class
$\cC_{\mu_i}$.  Any invariant connection $\ti{A}$ on $\ti{X} \times K$
descends to a connection $A$ on the quotient bundle $ (\ti{\mE}
\backslash \{ x_1,\ldots, x_b\} \times K) / \Gamma$ with holonomy
around $x_i$ in $\cC_{\mu_i}$.  Let
$\RR^\Gamma_K(\ti{X},\mu_1,\ldots,\mu_b)$ denote the moduli space of
$\Gamma$-equivariant flat bundles on $\ti{X}$ up to
$\Gamma$-equivariant isomorphism, such that the action of $\Gamma$ on
$\ti{\mE}_{x_{i}}$ is identified (up to conjugacy) with $\exp(\mu_i)$.
If $\alpha_0(\mu_i ) < 1 $ for $i = 1,\ldots, b$, the map
$(\ti{\mE},\ti{A}) \mapsto (\mE,A)$ defines a bijection
\begin{equation} \label{bundlethm}   
\RR^\Gamma_K(\ti{X},\mu_1,\ldots,\mu_b) \to \RR_K(X;\mu_1,\ldots,\mu_b) .\end{equation}
Indeed, any $\Gamma$-equivariant isomorphism of flat bundles on
$\ti{X}$ induces an isomorphism of bundles on $X \backslash \{ x_i
\}$.  Conversely, given a flat bundle on $X \backslash \{ x_i \}$ one
may pull-back to a flat bundle on $\ti{X} \backslash \{ \ti{x}_i \}$.
In polar coordinates $r_i,\theta_i$ near $x_i$ the connection has the
form $ \mu_i d \theta_i$.  It follows that one may glue in the trivial
flat bundle using the gluing map $\exp(\theta_i \mu_i)$ to obtain a
$\Gamma$-equivariant flat bundle on $\ti{X}$.  Any isomorphism of flat
bundles on $X \backslash \{ x_i \}$ lifts to an isomorphism of flat
bundles on $ \ti{X} \backslash \{ \ti{x}_i \}$.  In the local
trivializations near the marked points $x_i$ the isomorphism is given
by a constant gauge transformation $k$ in the centralizer of
$\exp(\mu_i)$, which is equal to the stabilizer of $\mu_i$ since
$\alpha_0(\mu_i) < 1$.  Therefore the isomorphism extends over the
points $\ti{x}_i$.

\subsection{The Yang-Mills heat flow}

According to Donaldson \cite{do:fo}, the Narasimhan-Seshadri
correspondence can be constructed by minimizing the Yang-Mills
functional on the space of connections.  Throughout this section we
identify the space $\A(X)$ of connections on $X \times K$ with the
space of holomorphic structures on $X \times G$.  For any connection
$A \in \A(X)$, let
$$d_A: \ \Omega^*(X,\k) \to \Omega^{* +1} (X,\k)$$
denote the corresponding covariant differentiation operator, and
$d_A^*$ its formal adjoint.  The Yang-Mills functional is $ A \mapsto
\Vert F_A \Vert^2_{L^2} .$ Let $\Theta$ be its contragradient, $
\Theta(A) = - d_A^* F_A .$ The connection $A$ is {\em Yang-Mills} if
$\Theta(A) = 0$.  The following summarizes results of Donaldson,
Daskalopolous, and R{\aa}de for $\Gl(n)$, extended to arbitrary
structure groups.

\begin{theorem} \label{ymthm}  
\begin{enumerate}
\item For any $A_0 \in \A(X)$, there exists a trajectory $A_t \in
C^0([0,\infty),\A(X))$ satisfying $ \ddt A_t = \Theta(A_t)$;
\item \label{limsec} $A_t$ converges in the Sobolev space $H^1$ to a
Yang-Mills connection $A_\infty$;
\item $A_\infty$ is a flat connection if and only if $A$ is
semistable; 
\item the map $A \mapsto A_\infty$ defines a continuous retract of the
space of semistable connections onto the space of flat connections, in
the $H^1$-topology; and
\item the map $[A] \mapsto [A_\infty]$ defines a homeomorphism
$\M_G(X) \to \RR_K(X)$.
\end{enumerate}
\end{theorem}

The results of \cite{do:fo,da:st,ra:ym} prove (a)--(e) for vector
bundles.  Fix an embedding $\phi: \ K \to U(n)$, and let the metric on
$\k$ be the pull-back of an invariant metric on $U(n)$.  The
Yang-Mills flow on $U(n)$-connections pulls back to the Yang-Mills
flow on $K$-connections.  This implies parts (a) and (b). Semistable
holomorphic structures on $X \times G$ map to semistable holomorphic
structures on $X \times \Gl(n)$, by functoriality of the canonical
reduction.  This implies (c) and (d).  It remains to show (e).

We must show that two connections $A_1,A_2$ are grade equivalent if
and only if the connections $A_{1,\infty}, A_{2,\infty}$ are in the
same $K(X)$-orbit.  First, we show that $S$-equivalence is equivalent
to topological $S$-equivalence, that is, $S$-equivalence where instead
of algebraic or holomorphic families of connections we require only
that the family be continuous, say in the Sobolev topology.  By
\cite[4.15.2]{ra:th}, there exists a (non-singular, projective)
universal space for semistable $G$-bundles on $X$, which we call
\label{UGsec} $\U_G(X)$ (Ramanathan's $R_3$.)  What we want to check
is that $\U_G(X)$ has the universal property for {\em topological}
families, at least locally.  That is, a continuous family $A_s$ of
semistable $G$-bundles defines a continuous family $B_s$ in $\U_G(X)$,
in a neighborhood of any $s_0 \in S$.  Choose an embedding $\iota: \ G
\to \Gl(V)$ and a line bundle $L \to X$, such that any bundle
$\iota_*(\mE) \otimes L$ is generated by globally sections and the
higher cohomology of $\iota_*(\mE) \otimes L$ vanishes.  A point in
$\U_G(X)$ is a set of generating sections for $\iota_*(\mE) \otimes L$,
together with a $G$-structure on $\iota_*(\mE)$.  Since higher
cohomology vanishes, the global sections of $\iota_*(A_s) \otimes L$
form a topological vector bundle over the parameter space $S$.  We can
choose a continuous family of sections $f_1(s),\ldots,f_m(s)$ that
generate $\iota_*(A_s)$ for any $s$ in a neighborhood $S_0$ of $s_0$.
Together with the $G$-structure on $\iota_*(A_s)$ these give the
family $B_s$.  Because $\M_G(X)$ is a good quotient of $\U_G(X)$, the
family $[A_s]$ is a continuous path in $\M$.  This shows that $A_0$
and $A_s$ are $S$-equivalent, for any $s \in S$.  In fact $\M_G(X)$ is
a coarse moduli space in the topological category, that is, represents
the functor from topological spaces to sets that assigns to any
topological space $S$ the set of continuous families $\mE_s, s \in S$ of
equivalence classes of semistable holomorphic $G$-bundles over $X$.
This implies that the holomorphic bundles corresponding to $A_j,
A_{j,\infty}$ are $S$-equivalent.  Hence, if $A_{1,\infty}$ are
isomorphic then $A_1,A_2$ are $S$-equivalent.

Conversely, suppose $A_1,A_2$ are semistable and $S$-equivalent.  The
grade bundles for $A_1,A_2$ are isomorphic, by \cite[3.12.1]{ra:th}.
Also, the grade bundles for $A_1,A_2$ are isomorphic to the grade
bundles of $A_{j,\infty}, \ j = 1,2$, since these bundles are
$S$-equivalent.  Since $A_{i,\infty}$ is flat, it is its own grade
bundle \cite[3.15]{ra:th}.  Flat connections isomorphic by a complex
gauge transformation are related by a unitary gauge transformation
\cite[Proposition 6.1.10]{do:fo}.  Hence $A_{1,\infty}$ and $
A_{2,\infty}$ are in the same $K(X)$-orbit, which completes the proof
of (e).

\subsection{Narasimhan-Seshadri theorems for equivariant
and parabolic bundles}

Let $\RR_{K,\Gamma}(\ti{X})$ denote the moduli space of
$\Gamma$-equivariant flat bundles on $\ti{X}$, up to equivariant
isomorphism.  Fix an action of $\Gamma$ on $\ti{X} \times K$.  If
$\ti{A}$ is a $\Gamma$-invariant connection on $\ti{X} \times K$, then
the tangent vector $\Theta(\ti{A})$ is also $\Gamma$-invariant.  The
Yang-Mills limit $\ti{A}_\infty$ is therefore a $\Gamma$-invariant
flat connection.  If $\ti{A}$ is semistable, then $\ti{A}_\infty$ is
flat, by \ref{ymthm} (c).  The map
\begin{equation} \label{eqcorresp}
\M_{G,\Gamma}(\ti{X}) \to \RR_{K,\Gamma}(\ti{X}) , \ \ [\ti{A}]
\mapsto [\ti{A}_\infty]\end{equation}
is a homeomorphism; the proof is essentially the same as in the
non-equivariant case.  This equivariant correspondence theorem implies
a correspondence theorem for parabolic bundles.  We need the following
lemma on existence of finite ramified covers.
\begin{lemma} \cite[5.2]{ed:re} \label{exists}
If $N$ is odd or $b$ is even, then there exists a $\Z_N$-cover $\pi: \
\ti{X} \to X$ totally ramified at $x_1,\ldots,x_b$.
\end{lemma}
Therefore, we can assume that $\ti{X}$ exists, at least after adding a
marked point with marking $0$.  

\begin{theorem} \label{corresp}  Let $G$ be a connected simple,
simply-connected Lie group with maximal compact subgroup $K$, and $X$
a curve with distinct marked points $x_1,\ldots,x_b$.  Let
$\mu_1,\ldots,\mu_b$ be markings with $\alpha_0(\mu_i) < 1$.  There is
a homeomorphism
$$ \M_G(X;x_1,\ldots,x_b,\mu_1,\ldots,\mu_b) \to \RR_K(X,\mu_1,
\ldots,\mu_b) .$$
\end{theorem}

For rational markings, this follows from Theorem \ref{isofunct} and
the bijections \eqref{eqcorresp} and \eqref{bundlethm}.  We extend it
to irrational markings by perturbation.  We note that Simpson's method
\cite{si:ha}, see also \cite{da:ymb} works just as well for the
non-rational case.

\begin{theorem}  \label{perturb} For any $(\mu_1,\ldots,\mu_b) \in \Alc^b$,
there exists a rational affine subspace $C(\mu_1,\ldots,\mu_b)$ such
that for $(\mu_1',\ldots,\mu_b')$ sufficiently close to
$(\mu_1,\ldots,\mu_b)$ in $C(\mu_1,\ldots,\mu_b)$ there exist
homeomorphisms
$$ \M_G(X;x_1,\ldots,x_b,\mu_1,\ldots,\mu_b) \to
\M_G(X;x_1,\ldots,x_b,\mu_1',\ldots,\mu_b') $$
\begin{equation} \label{rephom}
 \RR_K(X;\mu_1,\ldots,\mu_b) \to \RR_K(X;\mu_1', \ldots,\mu_b')
.\end{equation}
\end{theorem}

\begin{proof} 
Suppose $(\mu_1,\ldots,\mu_b)$ is contained in an open face
$F(\mu_1,\ldots,\mu_b) \subset \Alc^b$.  For each maximal parabolic
$P$, there exist a finite set 
$$S(\mu_1,\ldots,\mu_b) = \{ (d,w_1,\ldots,w_b) \in \N \times
(W/W_P)^b, \ \  d + \sum (w_i \mu_i, \omega_P) = 0 \} . $$
These inequalities define a rational affine subspace
$C(\mu_1,\ldots,\mu_b) \subset F(\mu_1,\ldots,\mu_b)$.  Let
$$m = \inf |d + \sum_{i=1}^b (w_i \mu_i, \omega_P)|, \ \ \ 
(d,w_1,\ldots,w_b) \notin S(\mu_1,\ldots,\mu_b) .$$ 
Since $W/W_P$ is finite and $\Alc$ is compact, $m$ is non-zero.  For
$(\mu_1',\ldots,\mu_b')$ sufficiently close to $(\mu_1,\ldots,\mu_b)$
in $C(\mu_1,\ldots,\mu_b)$ we have
$$ \M_G(X;x_1,\ldots,x_b,\mu_1,\ldots,\mu_b) =
\M_G(X;x_1,\ldots,x_b,\mu_1',\ldots,\mu_b') $$
since the semistability condition for the two sets of markings is the
same.

To prove the bijection for flat $K$-bundles, consider the manifold
$$M = K^{2(g + b- 1)} = \{
(a_1,\ldots,a_g,b_1,\ldots,b_g,c_1,\ldots,c_{b-1}, d_1,\ldots,d_{b-1})
\} $$
with action of $ (k_1,\ldots,k_b) \in K^b$ given by
$$ a_i \mapsto \Ad(k_b)a_i,\ b_i \mapsto \Ad(k_b)b_i,\ c_i \mapsto k_b
c_i k_i^{-1}, \ \ d_i \mapsto \Ad(k_i) d_i $$
and group valued moment map (see \cite{al:mom}; one could use loop
group actions here as well)
$$\Phi: \ M \to K^b, \ (a,b,c,d) \mapsto (d_1,\ldots,d_{b-1}, (\Pi
a_ib_ia_i^{-1} b_i^{-1} \Pi \Ad(c_i) d_i )^{-1}) .$$
The moduli space of flat bundles is the symplectic quotient
\begin{equation} \label{quot}
\RR_K(X;\mu_1,\ldots,\mu_b) = \Phinv(\mu_1,\ldots,\mu_b)/ (K_{\mu_1}
\times \ldots \times K_{\mu_b}) .\end{equation}
We claim that for $\nu \in TC(\mu_1,\ldots,\mu_b)$ and $\eps$
sufficiently small, there exists a homeomorphism
\begin{equation} 
 \RR_K(X;\mu_1,\ldots,\mu_b) \to \RR_K(X;\mu_1 + \eps \nu_1,
\ldots,\mu_b + \eps \nu_b) .\end{equation}
The quotient \eqref{quot} can be taken in stages, first a quotient by
$U(1)_\nu$ and then a quotient by $(K_{\mu_1} \times \ldots
K_{\mu_b})/U(1)_\nu$.  As in the Duistermaat-Heckman theorem, it
suffices to show that the one-parameter subgroup $U(1)_\nu$ generated
by $(\nu_1,\ldots,\nu_b)$ acts locally freely on
$\Phinv(\mu_1,\ldots,\mu_b)$.  Suppose $(a,b,c,d) \in
\Phinv(\mu_1,\ldots,\mu_b)$ is fixed by $U(1)_\nu$.  Then 
$$ a_i, b_i \in K_{\nu_b},\, d_i \in K_{\nu_i}, \, \nu_i = \Ad(c_i)
\nu_b .$$
Since $\nu_1,\ldots,\nu_b \in \t$, we have $\nu_i = \Ad(w_i) \nu_b$
for some $w_1,\ldots,w_b \in W$ and $c_i$ is a representative of
$w_i$, up to multiplication by $K_{\nu_b}$.  We may assume $\nu_b \in
\t_+$.  The fixed point set of $U(1)_\nu$ is
$$ K_{\nu_b}^{2g} \times (K_{\nu_b}w_1 \times K_{\nu_i}) \times \ldots
\times (K_{\nu_b} w_{b-1} \times K_{\nu_{b-1}}). $$
Its image under the moment map is equal to 
\begin{equation} \label{image}
 \{ (d_1,\ldots,d_b) \in K_{\nu_1} \times \ldots \times K_{\nu_b}, \ \
\ \prod_{i=1}^{b} \Ad(w_i) d_i \in [K_{\nu_b},K_{\nu_b}] \}
.\end{equation}
The stabilizer $K_{\nu_b}$ has roots $\alpha$ with $(\alpha,\nu_b) =
0$.  Therefore, the Cartan $\t_{\nu_b}$ of the semisimple part of
$K_{\nu_b}$ is
$$ \t_{\nu_b} = \on{span} \{ \alpha, \ (\alpha,\nu_b) =
0) \}  = \{ \xi \in \t, \ (\omega_j,\xi) = 0, \ j = 1,\ldots, m \}, $$ 
where $\omega_1,\ldots,\omega_m$ are the fundamental weights
corresponding to simple roots $\alpha_1,\ldots,\alpha_m$ with
$(\alpha_j,\nu_b) \neq 0$.  Let us identify the Weyl alcove $\Alc$
with a subset of $K$, using the exponential map.  The torus $T_{\nu_b}
\subset K_{\nu_b}$ intersects $\Alc$ in the subset defined by the
equations
$$ T_{\nu_b} \cap \Alc = \{ \xi \in \Alc, \ (\omega_j,\xi) \in \Z, \ j
= 1,\ldots,m \}.$$
The intersection of \eqref{image} with $\Alc^b$ is therefore
$$ \{ \xi \in \Alc^b, \ (\omega_j, \sum_{i=1}^b w_i \xi_i) \in \Z, \ \
j = 1,\ldots, m \} .$$
If $\xi = \mu $ belongs to this set then so does $\nu$, which implies
$$(\nu_b,\omega_j) = 0, \ j = 1,\ldots, m .$$  
Hence $\nu_b$ is a combination of simple roots vanishing on $\nu_b$
which is a contradiction.  \end{proof}

Working in the analytic category we can define a canonical
homeomorphism for non-rational markings
\begin{equation} \label{nonrat} \RR_K(X;\mu) \to \M_G(X;x;\mu) 
\end{equation} 
as follows.  Let $\ti{X} \to X$ denote the ramified cover with
covering group $\pi_1(X)$.  Any flat bundle on $X \backslash \{
x_1,\ldots, x_b \}$ with holonomies $\mu_1,\ldots,\mu_b$ defines a
$\pi_1(X)$-equivariant bundle on $\ti{X}$.  The corresponding
$\pi_1(X)$-equivariant holomorphic $G$-bundle $\ti{E} \to \ti{X}$
defines a parabolic $G$-bundle $E \to X$.  The resulting map
\eqref{nonrat} is continuous and injective, since the argument that
two flat bundles related by a complex gauge transformation are
unitarily isomorphic \cite[6.1.10]{do:fo} does not use compactness of
the curve.  Now consider the map
$$ \bigcup_{\nu \in U(\mu_1,\ldots,\mu_b)} \RR_K(X;\nu_1,\ldots,\nu_b)
 \to \bigcup_{\nu \in U(\mu_1,\ldots,\mu_b)}
 \M_G(X;x_1,\ldots,x_b;\nu_1,\ldots,\nu_b) $$
where $U(\mu_1,\ldots,\mu_b)$ is a closed neighborhood of
$(\mu_1,\ldots,\mu_b)$ in $C(\mu_1,\ldots,\mu_b)$ given by Theorem
\ref{perturb}.  Since this map is a homeomorphism for rational $\nu$
the image is dense.  In fact since the domain is compact, the map is
surjective.  It follows that the map is a homeomorphism, since the
domain is compact and the image is Hausdorff.  This completes the
proof of Theorem \ref{corresp}.

In this paper we do not deal with wall-crossing, that is, the change
in the topology of $\M(X;x;\mu)$ as the markings $\mu$ vary, see
\cite{do:va,th:fl} for the vector bundle case.

\section{Existence of parabolic bundles on the projective line}

We now turn to the question of which moduli spaces of parabolic
bundles are non-empty.  We continue to identify the space of
connections $\A(X)$ on $X \times K$ with the space of holomorphic
structures on $X \times G$.  Let $\A(X;x)$ denote the set of
holomorphic structures together with parabolic reductions at the
marked points $x_1,\ldots,x_b$, and let $\A(X;x;\mu)^{\ss}$ be the
subset corresponding to parabolic semistable bundles with markings
$\mu_1,\ldots,\mu_b$.  The moduli space $\M_G(X;x;\mu)$ is the
quotient of the $\A(X;x;\mu)^{\ss}$ by grade equivalence.  Let $f: \
\A(X;x) \to \A(X)$ denote the map forgetting the reductions.
$\A(X;x;\mu)^{\ss}$ is dense, if non-empty, in $\A(X;x)$.  For the
case without markings, this follows from Ramanathan's
\cite[5.8]{ra:th} or properties of the Shatz stratification
\cite{at:mo}.  The general case follows from the equivalence with
equivariant bundles.  Indeed, the equivalence shows that for any
finite-dimensional complex submanifold $S \subset \A(X;x;\mu)$,
$S^{\ss}$ is open and dense in $S$.  Any two points $A_1,A_2$ of
$\A(X;x;\mu)$ are contained in some $S$. By taking $A_1 \in
\A(X;x;\mu)^{\ss}$, one sees that $A_2 \in \ol{S^{\ss}} \subset
\ol{\A(X;x;\mu)^{\ss}}$.

\begin{lemma}  \label{gen}
For any markings $\mu_1,\ldots,\mu_b$, the moduli space $
\M_G(X;x;\mu)$ is non-empty if and only if the general element of
$\M_G(X;x;\mu)$ has a representative whose underlying principal bundle
is ordinary semistable.
\end{lemma}

\begin{proof} 
The intersection of a dense set with an open dense set is dense, hence
$\A(X;x;\mu)^{\ss} \cap \pi^{-1}(\A(X)^{\ss})$ is open and dense in
$\A(X;x;\mu)^{\ss}$.  Since the image of a dense set under a
surjective map is dense, the image of $ \A(X;x;\mu)^{\ss} \cap
\pi^{-1}(\A(X)^{\ss})$ is dense in $\M_G(X;x;\mu)$.
\end{proof}

We warn the reader that it is {\em not} true that the grade bundle of
a general element in the moduli space $\M_G(X;x;\mu)$ (that is, the
holomorphic bundle corresponding to a general element in
$\RR_K(X;\mu)$) is ordinary semistable.  For example, let $X = \P^1$,
$G = SL(2)$, and $\mu = (\mu_1,\mu_2,\mu_3)$ with $\mu_1 + \mu_2 +
\mu_3 = 1$.  Let $E = \P^1 \times \C$ be the trivial bundle with
general parabolic reductions at $x_1,x_2,x_3$.  Since there is a line
through three general points in $\P^1$, $E$ admits a parabolic
reduction with ordinary degree $1$ and parabolic degree $0$.  Hence
$\Gr(E)$ has underlying bundle $\O(1) \oplus \O(-1)$, which is
unstable.  The moduli space in this case is a single point, with
unitary representative $\Gr(E)$ which is not ordinary semistable.
Nevertheless, a general element of $\A(X;x)^{\ss}$ is ordinary
semistable.

\begin{proposition}  
\begin{enumerate}
\item
If $X$ has genus $g > 0$, then
$\M_G(X;x;\mu)$ is non-empty.
\item If $X$ has genus $ g= 0$, then
$\M_G(X;x;\mu)$ is non-empty if and only
if the trivial bundle $\mE = X \times G$ with general parabolic
structures at the marked points $x_1,\ldots,x_b$ is parabolic
semistable.
\end{enumerate}
\end{proposition}

\begin{proof}  
(a) follows from the holonomy description \eqref{holdes} and
surjectivity of the commutator map $K \times K \to K$ \cite{go:cm}.
(b) By Lemma \ref{gen}, $ \M_G(X;x;\mu)$ is
non-empty if and only if a semistable bundle $\mE$ with general
parabolic structures $\varphi_i$ is parabolic semistable.  By the
Birkhoff-Grothendieck theorem \cite{gr:cl} any principle $G$-bundle
admits a reduction of the structure group to $T_\C$.  Since $G$ is
simple, $c_1(\mE) = 0$.  A principal $T_\C$-bundle $\mE$ with $c_1(\mE) = 0$
is semistable if and only if $\mE$ is isomorphic to the trivial bundle,
which completes the proof.
\end{proof}

The trivial bundle $X \times G$ with parabolic structures
$(\varphi_i,\mu_i)$ is parabolic stable (resp. semistable) if and only
if
\begin{equation} \label{ineq2}
 \sum_{i=1}^b (w_i \omega_P, \mu_i) < d \ \text{ (resp. } \ \leq d
\text{ \ ) } \end{equation}
for all maximal parabolics $P$ and $([w_1],\ldots,[w_b]) \in W_i
\backslash W / W_P$ such that there exists a reduction $\sigma: \ X
\to G/P$ with $\deg(\sigma^* \mE(\omega_P)) = d$ with $\sigma(x_i)$ in
position $w_i$ relative to $\varphi_i$ for each $i = 1,\ldots,b$.

\begin{remark}  \label{add} If $(\mu_1,\ldots,\mu_b)$ is sufficiently
close to zero, then $\M_G(X;x;\mu)$ is isomorphic to the geometric
invariant theory quotient $(G/B)^b \qu G$.  More precisely, if
$$\sum_{i=1}^b (w_i \omega_P,\mu_i) < 1$$ 
for any maximal parabolic $P$ and any $(w_1,\ldots,w_b)$, then any
reduction $\sigma: \ X \to \mE/P$ with $\deg \sigma^* \mE(\omega_P) >
0$ violates semistability.  Therefore, in this case $\M_G(X;x;\mu)$
consists entirely of parabolic bundles whose underlying bundle
(forgetting the parabolic structure) is trivial.  Restricting the
condition \eqref{ineq2} to constant reductions $\sigma$ gives
precisely the stability condition for an element of $(G/B)^b$
\cite{be:coa}.  A symplectic argument for this fact is given in
Jeffrey \cite{je:em}.
\end{remark}

Since $\Del_b$ is a polytope of maximal dimension, it suffices to
consider the case that $\mu_1,\ldots,\mu_b$ are rational and lie in
the interior of $\Alc$.  In this case, $P_i = B$ for all $x_i$.
Suppose $ \varphi_i = g_i B $ for some $g_1,\ldots,g_b \in G$.  The
element $\sigma(x_i) \in G/P$ lies in position $w_i$ relative to
$\varphi_i \in G/B$ only if $\sigma(x_i)$ lies in the Schubert cell
$g_i C_{w_i}$.  Therefore,

\begin{proposition} The polytope $\Del_b$ is the set 
of points $(\mu_1,\ldots,\mu_b) \in \Alc^b$ satisfying
$$ \sum_{i=1}^b ( w_i \omega_P, \mu_i ) \leq d $$
for all maximal parabolics $P$ and $(w_1,\ldots,w_b) \in (W/W_P)^b$
such that there exists a holomorphic curve $\sigma: \ \P^1 \to G/P$
with $\sigma(x_i) \in g_iC_{w_i}$, for general $g_i \in G$.
\end{proposition}

We call an inequality \eqref{ineq2} essential if it actually defines a
facet (codimension one face) of $\Del_b$.  It remains to show that the
essential inequalities are those corresponding to the structure
coefficients $n_d(w_1,\ldots,w_b) = 1$.  The argument is the same as
that of Belkale \cite{bl:ip} in the vector bundle case.  Let
$(P;w_1,\ldots,w_b;d)$ define an essential inequality, and let
$(\mu_1,\ldots,\mu_b) \in \Alc^b$ be a point which violates that
inequality, and no others.  Let $\mE$ be a trivial $G$-bundle over
$\P^1$, with general parabolic structures $\varphi_i$ and markings
$\mu_i$.  Since $\mE$ is unstable, the canonical parabolic reduction
$\sigma_{\mE}$ is non-trivial, and defines an inequality which is violated
by $(\mu_1,\ldots,\mu_b)$.  Since only one inequality is violated,
$\sigma_{\mE}$ must be a reduction to a maximal parabolic, and the
corresponding inequality must be given by the data
$(P;w_1,\ldots,w_b;d)$.  The slope of $\sigma_{\mE}$ is $- d +
\sum_{i=1}^b (w_i \omega_P,\mu_i) .$ Since the canonical reduction is
the unique reduction of this slope, we must have $n_d(w_1,\ldots,w_b)
= 1 .$ This completes the proof of Theorem \ref{main}.

\section{The inequalities for type $G_2$}

We compute the small quantum cohomology for the generalized flag
varieties $G/P$ with $G$ of type $G_2$ and $P$ maximal, using a
variation on the quantum Chevalley formula of D. Peterson
\cite{pe:mi,fu:qp}.

\subsection{The quantum Chevalley formula}

Let $c_1(G/P) \in H^2(G/P) \cong \Z$ be the first Chern class of
$G/P$.  For any root $\beta$, let $h_\beta$ denote the corresponding
co-root.  Let $\alpha$ denote unique simple root, non-vanishing on
$\omega_P$.  

\begin{theorem}  For any $u \in W/W_P$ with minimal length representative
$\ti{u}$,
$$ [Y^{s_\alpha} ] \star [ Y^u] = \sum (h_\beta,\omega_P) [Y^{[\ti{u}s_\beta]}] +
\sum_{}
(h_\beta,\omega_P) q^{(h_\beta,\omega_P)} [Y^{[\ti{u}s_\beta]}] $$
where the first sum is over positive roots $\beta$ with
$l_P([\ti{u}s_\beta]) = l_P(u) + 1$ and the second is over positive
roots $\beta$ with $l_P([\ti{u}s_\beta]) = l_P(u) + 1 - c_1(G/P)
(h_\beta,\omega_P)$.
\end{theorem} 

\subsection{Small quantum cohomology for $G_2/P$, $P$ maximal}

Let $G$ be the group of type $G_2$, with simple roots
$\alpha_1,\alpha_2$ and positive roots
$$ \beta_1 = \alpha_1,\ \beta_2 = \alpha_2, \ \beta_3 = \alpha_1 +
\alpha_2,\ \beta_4 = 3 \alpha_1 + \alpha_2, \beta_5 \ = 2\alpha_1 +
\alpha_2, \ \beta_6 = 3\alpha_1 + \alpha_2 .$$
The highest root is $\beta_4$.  We fix the inner product on $\t^*$ so
that $(\beta_4,\beta_4) = 2$ and use it to identify $\t$ with $\t^*$.
The fundamental weights are $ \omega_1 = \beta_5, \ \ \omega_2 =
\beta_4 .$ The coroots and their pairings with the fundamental weights
are
$$ \begin{array}{c||c|c|} h_{\beta} & h_\beta(\omega_{1}) &
h_\beta(\omega_{2}) \\ \hline \hline 
h_{\beta_1} = 3\alpha_1  & 1 & 0 \\
h_{\beta_2} = \alpha_2 & 0 & 1 \\
h_{\beta_3} = 3\alpha_1 + 3\alpha_2 & 1 & 3 \\
h_{\beta_4} = 3\alpha_1 + 2 \alpha_2 & 1 & 2 \\
h_{\beta_5} = 6 \alpha_1 + 3 \alpha_2 & 2 & 3 \\
h_{\beta_6} = 3 \alpha_1 + \alpha_2 & 1 & 1 
\end{array} $$

Let $P_1,P_2$ denote the corresponding maximal parabolics, so that the
generalized flag varieties $G/P_j, j = 1,2$ have dimension $10$.  The
Weyl groups $W_{P_1},W_{P_2}$ are isomorphic to $\Z_2$.  Therefore,
the rational cohomology of $G/P_i$ is generated by a single generator
$y_1$ in degree $2$, with the single relation $y_1^6 = 0$.  The first
Chern classes in $H^2(G/P_i) \cong \Z$ are
$$ c_1(G/P_1) = 5, \ \ c_1(G/P_2) = 3 .$$
We denote by $y_i \in H_{10-2i}(G/P)$ the unique Schubert class of
codegree $2i$, with $y_0 = 1$.  The multiplication tables are given
below.  The second row in the table is given by Peterson's formula.
Since the cohomology is generated by $H^2$, the remaining rows in the
table may be computed recursively from the previous rows.

$$ \begin{array}{c||c|c|c|c|c|c|}
QH^*(G/P_1) & 1 & y_1 & y_2 & y_3 & y_4 & y_5 \\
\hline \hline
1 & 1 & y_1 & y_2  & y_3 & y_4     & y_5 \\
y_1 &     &y_2 & 2y_3 & y_4 & y_5 + q & qy_1 \\
y_2 &     &    & 2y_4 & y_5 + q & 2qy_1 & qy_2 \\
y_3 &     &    &      & qy_1 & qy_2 & qy_3 \\
y_4 &     &    &      &      & 2qy_3& qy_4 \\
y_5 &     &    &      &      &      & q^2 
\end{array}
$$

$$ \begin{array}{c||c|c|c|c|c|c|}
QH^*(G/P_2) & 1 & y_1 & y_2 & y_3 & y_4 & y_5 \\
\hline \hline
1 & 1 & y_1 & y_2        & y_3          & y_4            & y_5 \\
y_1 &     & 3y_2&2y_3+q      &3y_4+qy_1     & y_5+qy_2       & qy_3 + 2q^2 \\
y_2 &     &     & 2y_4+ qy_1 & y_5 + 2qy_2  & qy_3 + q^2     & qy_4 + q^2y_1\\
y_3 &     &     &            & 2qy_3 + 2q^2 & qy_4 + q^2y_1  & 2q^2 y_2\\
y_4 &     &     &            &              & q^2 y_2        & q^2 y_3  \\
y_5 &     &     &            &              &                & 2q^2 y_4        \end{array}
$$

From the second row of the tables one may also compute the presentation
of the quantum cohomology rings, in terms of the generator $y = y_1$.
First, one obtains the following Giambelli-type expressions for the
Schubert classes.  In $QH^*(G/P_1)$, we have
$$ y_2 = y_1^2, \ y_3 = y_1^3/2, \ y_4 = y_1^4/2, \ y_5 = y_1^5/2 - q $$
and in $QH^*(G/P_2)$
$$ y_2 = y_1^2/3, \ y_3 = (y_1^3 - 3q)/6, \ y_4 = (y_1^4 - 9qy_1)/18, \ 
y_5 = (y_1^5 - 15 qy_1^2)/18 .$$

From the last entry in the first row in the tables one obtains
\begin{proposition}  $QH^*(G/P_1)$ is generated by a single generator
$y$ of degree $2$, with relation $y_1^6 = 4qy_1$.  $QH^*(G/P_2)$ is
generated by a single generator $y$ of degree $2$, with relation $y_1^6
= 18qy_1^3 + 9q^2 $.
\end{proposition}

In particular, both of these rings are semisimple at $q=1$, since the
relations have no multiple roots.  Neither the classical integral or
quantum rational cohomology of $G/P_1$ and $G/P_2$ is the same as that
of complex projective space $\C P^5$, although the classical rational
cohomology is the same.

\subsection{The inequalities}

From the tables, one may read off 33 classical and 40 quantum
inequalities.  Some of the quantum inequalities do not define facets;
it would be interesting to determine which ones.  For example, the
last entry in the table for $G/P_1$ gives the inequality $ ( \omega_1,
\mu_1 + \mu_2 + \mu_3) \leq 2 $ which does not define a facet since
$(\mu,\omega_1) \leq (\omega_1,\omega_1) = 2/3$ for any $\mu \in
\Alc$.

\def\cprime{$'$} \def\cprime{$'$} \def\cprime{$'$}
  \def\polhk#1{\setbox0=\hbox{#1}{\ooalign{\hidewidth
  \lower1.5ex\hbox{`}\hidewidth\crcr\unhbox0}}}

\date{Revised January 6, 2003}

\end{document}